\documentclass[11pt,amsfonts]{amsart}

\usepackage{fullpage}

\begin{document}

\newtheorem{theorem}{Theorem}[section]
\newtheorem{corollary}[theorem]{Corollary}
\newtheorem{definition}[theorem]{Definition}
\newtheorem{lemma}[theorem]{Lemma}
\newtheorem{proposition}[theorem]{Proposition}
\newtheorem{remark}[theorem]{Remark}

\def\endproof{\qed\medskip}
\def\blacksquare{\hbox to .60em {\vrule width .60em height .60em}}
\renewcommand{\qedsymbol}{\hfill \blacksquare \hspace *{4.5em}}

\renewcommand{\theequation}{\thesection.\arabic{equation}}

\title[]{Dehn Filling and Einstein Metrics in Higher Dimensions}

\author[]{Michael T. Anderson}

\thanks{$^{*}$ Partially supported by NSF Grants DMS 0072591 and 0305865}

\maketitle

\abstract
We prove that many features of Thurston's Dehn surgery theory for 
hyperbolic 3-manifolds generalize to Einstein metrics in any dimension. 
In particular, this gives large, infinite families of new Einstein 
metrics on compact manifolds.
\endabstract

\setcounter{section}{0}
\setcounter{equation}{0}

\section{Introduction.}

  In this paper, we construct a large new class of Einstein metrics of 
negative scalar curvature on $n$-dimensional manifolds $M = M^{n}$, for 
any $n \geq 4$. Einstein metrics are Riemannian metrics $g$ of constant 
Ricci curvature, and we will assume the curvature is normalized as 
\begin{equation}\label{e1.1}
Ric_{g} = -(n-1)g,
\end{equation}
so that the scalar curvature $s = -n(n-1)$.
  The construction is a direct generalization of Thurston's theory of 
Dehn surgery or Dehn filling on hyperbolic 3-manifolds [31] to Einstein 
metrics in any dimension; in fact the proof gives a new, analytic 
approach to Thurston's cusp closing theorem [31], [32].

\medskip

  To describe the construction, start with any complete, non-compact 
hyperbolic $n$-manifold $N = N^{n}$ of finite volume, with metric 
$g_{-1}$ of constant curvature $-1$. The manifold $N$ has a finite 
number of cusp ends $\{E_{j}\}$, $1 \leq j \leq q$, with each end $E$ 
diffeomorphic to $F\times {\mathbb R}^{+}$, where $F$ is a compact flat 
manifold, with flat metric $g_{0}$ induced from $(N, g_{-1})$. For 
simplicity, assume that each $F$ is an $(n-1)$-torus $T^{n-1}$; this 
can always be achieved by passing to a finite covering space if 
necessary, cf. [6].

  Now perform Dehn filling on any collection ${\mathcal C} = \{E_{k}\}$ 
of cusp ends of $N$, where $1 \leq k \leq p$, and $p \leq q$. Thus, 
fix a torus $T^{n-1} \subset E \in {\mathcal C}$ and let $\sigma$ be a 
simple closed geodesic $\sigma \subset (T^{n-1}, g_{0})$. Attach a 
(generalized) solid torus $D^{2}\times T^{n-2}$ onto $T^{n-1}$ by a 
diffeomorphism of $\partial D^{2}\times T^{n-2} \simeq T^{n-1}$ sending 
$S^{1} = \partial D^{2}$ onto $\sigma$. If $\sigma_{k}$ are such simple 
closed geodesics in tori $T_{k}^{n-1} \subset E_{k}$, let ${\bar \sigma} 
= (\sigma_{1}, ... \sigma_{p})$ and let 
\begin{equation} \label{e1.2}
M = M_{\bar \sigma} = M^{n}(\sigma_{1}, ..., \sigma_{p})
\end{equation}
be the resulting manifold obtained by Dehn filling the collection of ends 
$E_{1}, ..., E_{p}$ of $N$. The diffeomorphism type of $M$ depends on 
the homotopy class of each $\sigma_{k}$ in $\pi_{1}(T_{k}^{n-1}) 
\simeq {\mathbb Z}^{n-1}$ but is otherwise independent of the choice 
of attaching map.

  If $p = q$ the manifold $M_{\bar \sigma}$ is compact, (without boundary); 
otherwise $M_{\bar \sigma}$ has $q - p$ remaining cusp ends. Define the 
Dehn filling ${\bar \sigma} = (\sigma_{1}, ..., \sigma_{p})$ 
to be {\it sufficiently large} if, given $N$ and a fixed collection of 
tori $T_{k}^{n-1}$, the length $R_{k}$ of each geodesic $\sigma_{k}$, 
$1 \leq k \leq p$, is sufficiently large in $(T_{k}^{n-1}, g_{0})$; 
this will be made more precise in \S 3.  

  The main result of the paper is then the following:

\medskip

\begin{theorem}\label{t1.1}
 Let $(N, g_{-1})$ be a complete, non-compact hyperbolic $n$-manifold 
of finite volume, $n \geq 3$, with toral ends. Then any manifold 
$M_{\bar \sigma}$ obtained by a sufficiently large Dehn filling 
of the ends of $N$ admits a complete, finite volume  Einstein metric $g$, 
of uniformly bounded curvature and satisfying (1.1).
\end{theorem}

  To place this result in some perspective, a well-known result of Wang 
[33] states that if $n \geq 4$, there are only finitely many complete 
hyperbolic $n$-manifolds with volume $\leq V$. On the other hand, let 
${\mathcal H}(V)$ denote the number of (diffeomorphically) distinct 
complete non-compact hyperbolic $n$-manifolds of volume $\leq V$. Then 
${\mathcal H}(V)$ grows super-exponentially in $V$; in fact, by a 
recent result in [13], there are constants $a$, and $b$, depending only 
on $n$, such that
\begin{equation}\label{e1.3}
e^{aV\ln V} \leq {\mathcal H}(V) \leq e^{bV\ln V}.
\end{equation}
(The lower bound in (1.3) is stated in [13] only for compact hyperbolic 
manifolds, but using the work of Lubotzky in [24], this bound also 
holds for non-compact hyperbolic manifolds, [25]). For many such 
manifolds $N$, the number of cusp ends also grows linearly in $V$, cf. 
Remark 4.2.

  With each such $N$, Theorem 1.1 associates infinitely many 
homeoeomorphism types of compact manifolds $M_{\bar \sigma}$, 
(as well as non-compact manifolds). Formally, the number 
of such compact manifolds is $\infty^{q}$, where $q$ is the number of cusps 
of $N$. The Einstein metrics all have volume close to $V = vol N$. Further, 
although all hyperbolic manifolds are locally isometric, most of the Einstein 
metrics constructed are not locally isometric. Thus, the result gives a wealth 
of new examples of Einstein manifolds. 

\medskip

  All of the manifolds $M_{\bar \sigma}$ are $K(\pi, 1)$ manifolds, 
again for $\bar \sigma$ sufficiently large; in fact all admit metrics 
of non-positive sectional curvature. However, none of these manifolds 
admit metrics of negative sectional curvature. The curvature of the 
Einstein metrics $g$ on $M_{\bar \sigma}$ is not non-positive, (at least when 
$n > 4$), but one has the uniform bounds
\begin{equation}\label{e1.4}
-1 - \tfrac{1}{2}(n-3) - \varepsilon({\bar \sigma}) \leq K \leq -1 + 
\tfrac{1}{2}(n-3)(n-2) + \varepsilon({\bar \sigma}),
\end{equation}
where $K$ denotes the sectional curvature of the metric, and 
$\varepsilon({\bar \sigma})$ is small, with $\varepsilon({\bar \sigma}) 
\rightarrow 0$ as ${\bar \sigma} \rightarrow \infty$ in the Dehn 
filling space attached to each cusp. When $n = 4$, note that (1.4) gives 
$K \leq \varepsilon(\bar \sigma)$, so that the Einstein metrics are of almost 
non-positive curvature. When $n = 3$, the Einstein metrics are of course hyperbolic; 
the construction in Theorem 1.1 then gives an analytic proof of Thurston's 
cusp closing theorem. 

\medskip

  The Einstein metrics $(M, g)$ given by Theorem 1.1 are all 
close to the initial hyperbolic manifold $(N, g_{-1})$ in the pointed 
Gromov-Hausdorff topology. This will be apparent in a precise sense 
from their construction, but can be formulated generally as follows. 
Note first that $N$ is embedded in any $M$ obtained by Dehn 
filling as the complement of a generalized link - the collection of the 
$(n-2)$ tori at the core of the solid tori $\{D_{j}^{2}\times T_{j}^{n-2}\} 
\subset M$. Given $(N, g_{-1})$, let $g^{k}$ be a sequence of Einstein metrics 
on $M^{k} = M_{\bar \sigma^{k}}$, (constructed by the Theorem), such that the 
length of $\sigma_{j}^{k}$ diverges to infinity as $k \rightarrow 
\infty$, for each $\sigma_{j}^{k} \in \bar \sigma^{k}$. Then, given a 
fixed base point $y \in N \subset M^{k}$, the metrics $(M^{k}, g^{k}, 
y)$ converge to $(N, g_{-1})$ in the pointed Gromov-Hausdorff topology 
based at $y$. 

  The convergence is smooth on compact domains containing $y$, and the 
curvature tends to $-1$, uniformly on compact subsets. Thus, by the 
bounds (1.4) and the fact that the volume of $(M^{k}, g^{k})$ is 
uniformly bounded, one finds that the metrics $(M^{k}, g^{k})$ have 
uniformly small Weyl curvature in $L^{p}$, for any $p < \infty$: 
\begin{equation}\label{e1.5}
\int_{M^{k}}|W|^{p}dV_{g^{k}} \leq \varepsilon({\bar \sigma^{k}}, p),
\end{equation}
where $\varepsilon$ depends only on $p$ and $\bar \sigma^{k}$; for any 
fixed $p$, $\varepsilon \rightarrow 0$ as the length of 
$\sigma_{j}^{k}$ diverges to infinity, for all $j$. This behavior does 
not hold w.r.t. the $L^{\infty}$ norm. 

  We also point out that each Einstein metric $g$ constructed on any $M 
= M_{\bar \sigma}$ is an isolated point in the moduli space of Einstein 
metrics on $M$, cf. Remark 3.8; thus such metrics are (locally) rigid.

\medskip

  Theorem 1.1 is an analogue of Thurston's cusp closing theorem [31]. The 
next result is an analogue of the Jorgensen-Thurston cusp opening theorem, cf. 
[31], [18]. Let ${\mathcal E}$ be the class of complete, finite volume Einstein 
metrics constructed via Theorem 1.1, together with the class of complete, non-compact 
hyperbolic $n$-manifolds $(N, g_{-1})$ of finite volume. Let ${\mathcal E}_{V}$ 
be the subset of ${\mathcal E}$ of metrics of volume $\leq V$.
\begin{theorem}\label{t1.2}
The space ${\mathcal E}$ is closed with respect to the pointed Gromov-Hausdorff 
and $C^{\infty}$ topologies and the subspaces ${\mathcal E}_{V}$ are compact, 
for any $V < \infty$. Any limit point $(M_{\infty}, g_{\infty}) \in {\mathcal E}$ 
of a sequence $(M^{k}, g^{k}) \in {\mathcal E}$ satisfies
$$C(M_{\infty}) > \max_{k} C(M^{k}),$$
where $C(M)$ denotes the number of cusp ends of $M$. 
\end{theorem}

  In fact, Theorem 1.1 is proved for compact manifolds, where one Dehn-fills all 
the cusp ends of a given hyperbolic manifold $N$. It is then shown that the closure 
of the class of resulting Einstein metrics in the pointed Gromov-Hausdorff topology 
consists of the Einstein manifolds satisfying the conclusions of Theorem 1.1. Given 
this, the main content of Theorem 1.2 is the compactness ${\mathcal E}_{V}$. 

  Taken together, these results are close analogues of Thurston's Dehn 
surgery theory for hyperbolic 3-manifolds. Note that the manifolds $M_{\bar \sigma}$ 
in Theorem 1.1 can be viewed as obtained by Dehn {\it surgery} on a fixed manifold 
$M = M_{\bar \sigma_{0}}$, where $\bar \sigma_{0}$ is any Dehn filling of 
all the ends of $N$. The original non-compact hyperbolic manifold $N$ is 
then given by $N = M(\infty, ..., \infty)$.

\medskip

  Several aspects of the Thurston-Jorgensen picture of the structure of 
the volumes of hyperbolic 3-manifolds also generalize to Einstein metrics 
in higher dimensions. We describe briefly here the picture in dimension 4; 
further details, and discussion of the volume behavior in higher 
dimensions, are given in \S 4.

  The Chern-Gauss-Bonnet theorem shows that the volume of a complete, 
finite volume hyperbolic 4-manifold is given by
\begin{equation}\label{e1.6}
vol(N, g_{-1}) = \frac{4\pi^{2}}{3}\chi(N) \geq 0.
\end{equation}
Further, it is known that given any $k \in {\mathbb Z}^{+}$, there are 
(many) complete, non-compact hyperbolic $4$-manifolds $N^{k}$ of finite 
volume, with $\chi(N^{k}) = k$, cf. [29] for example. Let $(M, g)$ be 
any Einstein metric constructed via Theorem 1.1. Then the 
Chern-Gauss-Bonnet theorem gives
\begin{equation}\label{e1.7}
vol(M, g) = \frac{4\pi^{2}}{3}\chi(M) - \frac{1}{6}\int_{M}|W|^{2}.
\end{equation}
By a standard Mayer-Vietoris argument, $\chi(M) = \chi(N)$ and thus by 
(1.6),
\begin{equation}\label{e1.8}
vol(M, g) = vol(N, g_{-1}) - \delta(\bar \sigma) < vol(N, g_{-1});
\end{equation}
here $\delta(\bar \sigma)$ is small, and by (1.5) may be made 
arbitrarily small if the Dehn fillings in $\bar \sigma = 
(\sigma_{1}, ..., \sigma_{p})$ are sufficiently large, depending on 
$\delta$. Thus, the volume decreases under Dehn filling. 

  Several features of the Thurston-Jorgensen theory of volumes of hyperbolic 
3-manifolds thus generalize to Einstein metrics in higher dimensions. In 
particular, the set of volumes of metrics in ${\mathcal E}$ is a non-discrete, 
countable closed set in ${\mathbb R}$. However, it is not known if the set of 
volumes is well-ordered, (as a subset of $\mathbb R$, or finite-to-one, as in 
the Thurston-Jorgensen theory; again see \S 4 for further discussion.

\medskip

  The main idea of the proof is a glueing procedure now frequently used 
in constructing solutions to geometric PDE. Thus, one constructs an 
approximate Einstein metric on $M = M_{\bar \sigma}$, and shows this 
can be perturbed to an exact solution, i.e. an Einstein metric, by means 
of the inverse function theorem. Most of the technical work in 
the paper is concerned with the proof that the linearization of the 
Einstein operator (1.1) uniformly near the approximate solution is an 
isomorphism, (modulo diffeomorphisms). 

  Conceptually, the main issue is to construct the approximate 
solution. Since the hyperbolic manifold $N$ is already Einstein, one 
needs to find suitable complete Einstein metrics on $D^{2}\times 
T^{n-2}$ which asymptotically approach a hyperbolic cusp metric. Now a 
model for such metrics was constructed long ago by physicists, see [23] 
for instance, and later by Berard-Bergery [7], cf. also [8,9.118]. More 
recently these model metrics have been frequently analysed in 
connection with the AdS/CFT correspondence, and are now commonly called 
toral AdS black hole metrics, cf. [12] and references therein for 
example. These metrics have the following simple explicit form:
\begin{equation}\label{e1.9}
g_{BH} = V^{-1}dr^{2} + Vd\theta^{2} + r^{2}g_{T^{n-2}},
\end{equation}
where $g_{T^{n-2}}$ is any flat metric on $T^{n-2}$ and $V = V_{m}(r)$ 
is the function
\begin{equation}\label{e1.10}
V = r^{2} - \frac{2m}{r^{n-3}},
\end{equation}
If $n = 3$, this gives the usual hyperbolic metric on a tube about a 
single core geodesic. The parameter $r$ runs over the interval $[r_{+}, 
\infty)$, where $r_{+} = (2m)^{1/n-1}$. In order to obtain a smooth 
metric, the circular parameter $\theta$ is required to run over the 
interval $[0, \beta]$, where
\begin{equation}\label{e1.11}
\beta = \frac{4\pi}{(n-1)r_{+}}.
\end{equation}
The number $m$ is any positive number, and represents the mass of 
$g_{BH}$.

  The metric $g_{BH}$ has infinite volume, and so is not asymptotic to 
a hyperbolic cusp in the usual sense. However, we will see that this 
can be remedied by suitably ``twisting'' these metrics. This has been 
previously described in [1] and is discussed further in \S 2 below. 
Briefly, all the metrics $g_{BH}$ in (1.9) are isometric in the 
universal cover $D^{2}\times {\mathbb R}^{n-2}$. By taking suitable 
isometric actions of ${\mathbb Z}^{n-2}$ on the universal cover, the 
quotient has large regions closely approximating a given hyperbolic 
cusp metric. Thus, one may glue on a suitable quotient of the metric 
$g_{BH}$ onto a cusp of $N$ to obtain an approximate Einstein metric. 
This is exactly the same observation as Thurston's in the context of 
Dehn filling of hyperbolic 3-manifolds.

\medskip

  There is a large and growing literature on such glueing constructions 
for numerous geometric PDE. However, these have not been previously 
successful in constructing Einstein metrics; to our knowledge, the only 
exception is the work of Joyce on the construction of Einstein metrics  
of special holonomy in dimensions 7 and 8. More recently, Mazzeo and 
Pacard [26] have constructed new classes of conformally compact Einstein 
metrics on open manifolds, (of infinite volume), by a glueing technique 
on the boundary at conformal infinity.

\medskip

  The contents of the paper are briefly as follows. In \S 2, we discuss 
a number of background results and material needed for the proof of 
Theorem 1.1. The proof of Theorem 1.1 follows in \S 3. Several further 
results are then given after the proof. Thus, Proposition 3.9 proves 
that there are only finitely many Dehn fillings of a given $N$ which 
have the same homeomorphism type, while Corollary 3.11 discusses Dehn 
fillings on non-toral ends. In \S 4, we discuss a number of aspects of 
the geometry and topology of the manifolds $M_{\bar \sigma}$, as well 
as the convergence and volume behavior of the set of all Einstein 
metrics constructed by Dehn filling. Theorem 1.2 is proved at 
the end of \S 4.1.

\medskip

  I would like to thank Lowell Jones, Alex Lubotzky and Pedro Ontaneda 
for their assistance and discussions on hyperbolic manifolds and homotopy 
equivalences, and Gordon Craig for his assistance with the manuscript. 
Thanks also to Claude LeBrun and Dennis Sullivan for their comments and 
interest in this work. I especially thank Rafe Mazzeo for enlightening 
discussions on the behavior of Einstein metrics on cusp-like ends 
and for insightful comments and criticism of various aspects of the 
paper. Finally, my thanks to the referees for their careful examination 
and very constructive comments on the paper.

\section{Background Material.}

\setcounter{equation}{0}

\medskip

  In this section, we assemble background results and material needed 
for the work in \S 3. We break the discussion into four subsections 
dealing with different topics. 

\medskip

{\bf \S 2.1.}  Let $(N, g_{-1})$ be a complete, non-compact hyperbolic 
manifold of finite volume. As mentioned in the Introduction, $N$ then 
has a finite number of ends $E_{i}$, each diffeomorphic to $F\times 
{\mathbb R}^{+}$, where $F$ is a flat manifold; the topological type of 
$F$ depends of course on the end $E$. 

  It is not difficult to show that there is a finite cover $\bar N$ of 
$N$ such that all ends of $\bar N$ are tori $T^{n-1}$, cf. [6, Cor. 
2.4] for instance. For simplicity, from now on, we assume this is the 
case, and drop the bar from the notation; see Lemma 3.10 for discussion 
of non-toral ends. 

  The groups $\pi_{1}(T^{n-1}) \simeq {\mathbb Z}^{n-1}$ inject in 
$\pi_{1}(N)$ and are called the peripheral subgroups of $\pi_{1}(N)$. 
Any subgroup of $\pi_{1}(N)$ isomorphic to ${\mathbb Z}^{n-1}$ is 
conjugate to some peripheral subgroup; in fact any non-cyclic abelian 
subgroup is conjugate to a subgroup of some peripheral subgroup. 

\medskip

   The hyperbolic metric $g_{-1}$ on any cusp end $E$ has the form
\begin{equation}\label{e2.1}
g_{-1} = dt^{2} + e^{2t}g_{0},
\end{equation}
where $g_{0}$ is a flat metric on the $(n-1)$-torus $T^{n-1}$, and $t$ 
runs over the interval $(-\infty, 0]$. By the Margulis Lemma [18], 
[22], the flat metric $g_{0}$ may be chosen so that the injectivity 
radius $inj_{g_{0}}$ satisfies $inj_{g_{0}}T^{n-1} \geq \mu_{0}$, for a 
fixed constant $\mu_{0}$, depending only on $n$. For each end $E$ of 
$N$ on which Dehn filling is performed, we thus choose a fixed toral 
slice $T^{n-1} = \{0\} \times T^{n-1} \subset E$ satisfying this 
property. Given this, one may then write
$$(T^{n-1}, g_{0}) = {\mathbb R}^{n-1} / {\mathbb Z}^{n-1},$$
where the lattice ${\mathbb Z}^{n-1}$ is generated by $(n-1)$ basis 
vectors $v_{1}, ..., v_{n-1} \in {\mathbb R}^{n-1}$. The vectors 
$v_{i}$ are naturally identified with simple closed geodesics in 
$(T^{n-1}, g_{0})$ which intersect each other exactly once in a single 
base point. The choice of lattice vectors $(v_{1}, ..., v_{n-1})$ is of 
course not unique - it may be changed by any element in $SL(n-1, 
{\mathbb Z})$. However, we again fix such a basis of each 
$\pi_{1}(T^{n-1})$ once and for all.

\medskip

  Next, we describe the process of Dehn filling in higher dimensions; 
this is completely analogous to the situation in 3 dimensions. 

  Fix an end $E$ and $T^{n-1} \subset E$ as above. Elements $[\sigma]$ 
of $\pi_{1}(T^{n-1}) \simeq {\mathbb Z}^{n-1}$ are represented by 
closed geodesics in $(T^{n-1}, g_{0})$. If $\sigma$ is then any simple 
closed geodesic in $(T^{n-1}, g_{0})$, the class $[\sigma]$ may be 
represented in the form
$$[\sigma] = \sum \sigma^{i}[v_{i}],$$
where each $\sigma^{i} \in {\mathbb Z}$ and the collection $\sigma^{I} 
= (\sigma^{1}, ..., \sigma^{n-1})$ is primitive, in the sense that 
$\sigma^{I}$ is not a multiple of some $\sigma^{I'}$. 

  Now attach a (generalized) solid torus $D^{2}\times T^{n-2}$ to 
$T^{n-1}$ by a diffeomorphism $\phi$ of the boundary $\partial 
(D^{2}\times T^{n-2}) = S^{1}\times T^{n-2}$ with $T^{n-1}$, which 
sends $S^{1}$ to the closed geodesic $\sigma$. This gives the Dehn 
filled manifold
\begin{equation}\label{e2.2}
M_{\sigma} = (D^{2} \times T^{n-2})\cup_{\phi} N.
\end{equation}
By the Bieberbach rigidity theorem [9], any diffeomorphism of $T^{n-1}$ 
is isotopic to an element of $SL(n-1, {\mathbb Z})$, and so extends to 
a diffeomorphism of the solid torus $D^{2}\times T^{n-2}$. Thus the 
topological type of $M_{\sigma}$ is well-defined by the homotopy class 
of $[\sigma] \in \pi_{1}(T^{n-1})$. In fact, the topological type of 
$M_{\sigma}$ depends only on the unoriented curve $\sigma$, i.e. the 
class $[\pm \sigma] \in \pi_{1}(T^{n-1})$, cf. [30]. The vector 
$$ \sigma = (\sigma^{1}, ..., \sigma^{n-1})$$
gives the filling coefficients associated to $\sigma$, (w.r.t. the 
basis $\{v_{i}\}$). The Dehn filling space associated to the end $E$ is 
the collection of primitive $(n-1)$-tuples $\{\sigma^{i}\}$, and thus a 
subset of ${\mathbb Z}^{n-1}/ \{\pm 1\}$.

  This process may be carried out separately on any collection of ends 
$E_{j}$, $1 \leq j \leq p \leq q$, of $N$ and gives the manifold $M = 
M_{\bar \sigma}$, ${\bar \sigma} = (\sigma_{1}, ..., \sigma_{p})$, 
obtained by Dehn filling on the ends of $N$.

\medskip

  Next we make a number of remarks on the topology of the manifolds $M 
= M_{\bar \sigma}$. First, the hyperbolic manifold $N$ embeds in any 
$M$,
\begin{equation}\label{e2.3}
N \subset M
\end{equation}
as the complement of the core tori $T^{n-2}$ of each Dehn filling. We 
recall a well-known result of Gromov-Thurston, the $2\pi$ theorem, cf. 
[19] or [11,Thm.7]; this states that when the length $L(\sigma)$ of 
$\sigma$ in the flat torus $(T^{n-1}, g_{0})$ satisfies
\begin{equation}\label{e2.4}
L(\sigma) \geq 2\pi,
\end{equation}
the resulting manifold $M_{\sigma}$ has a complete metric of 
non-positive sectional curvature and finite volume. Although proved in 
the context of 3-manifolds, the same result and proof holds in any 
dimension. (Briefly, one forms the Euclidean cone of length 1 on 
$\sigma$, and takes the constant skew product with the flat metric on 
$T^{n-2}$. This gives a singular flat metric on $D^{2}\times T^{n-2}$, 
with cone angle $L(\sigma)$ along the core $T^{n-2}$. A natural 
smoothing of this cone singularity gives a metric of non-positive 
curvature on $M_{\sigma}$). 

  In particular, all the manifolds $M_{\bar \sigma}$ satisfying (2.4) 
for each geodesic $\sigma_{j} \in \bar \sigma$ are $K(\pi, 1)$ 
manifolds. Further, with respect to the metric of non-positive 
curvature on $M_{\bar \sigma}$, the core tori $T^{n-2}$ are totally 
geodesic. Since all closed geodesics in a manifold of non-positive 
curvature are essential in $\pi_{1}$, it follows that each core torus 
injects in $\pi_{1}$:
\begin{equation}\label{e2.5}
\pi_{1}(T^{n-2}) \hookrightarrow \pi_{1}(M_{\bar \sigma}).
\end{equation}
In particular, by Preissman's theorem, one sees that $M_{\bar \sigma}$ 
does not admit a metric of negative sectional curvature when dim $M 
\geq 4$.

\bigskip

{\bf \S 2.2.}  In this subsection, we discuss some aspects of the 
geometry of the (standard) AdS toral black hole metrics (1.9):
\begin{equation}\label{e2.6}
g_{BH} = V^{-1}dr^{2} + Vd\theta^{2} + r^{2}g_{T^{n-2}}.
\end{equation}
As in (1.10) and (1.11), $V = V(r) = r^{2} - 2mr^{-(n-3)}$ and $\theta$ 
takes values in $[0, \beta]$, where $\beta = 4\pi / (n-1)r_{+}$, $r_{+} 
= (2m)^{1/(n-1)}$ with $r \in [r_{+}, \infty)$. Although this metric 
appears to be singular at $r = r_{+}$, a simple change of coordinates, 
(analogous to the change from polar to Cartesian coordinates), shows 
that $g_{BH}$ is smooth everywhere. The metric is defined on the solid 
torus $D^{2} \times T^{n-2}$ and $g_{T^{n-2}}$ is any flat metric on 
$T^{n-2}$. 

  From the physical point of view, the core $(n-2)$-torus $H = \{r = 
r_{+}\} \subset D^{2}\times T^{n-2}$ represents the horizon of a black 
hole. Note that $H$ is the fixed point set of the isometric $S^{1}$ 
action given by rotation in $\theta$. Thus, $H$ is totally geodesic in 
$g_{BH}$; $H$ gives the usual core geodesic in a hyperbolic tube when 
$n = 3$. 

  The metric $g_{BH}$ is an Einstein metric, satisfying (1.1), which is 
asymptotically hyperbolic or conformally compact, cf. [2] or [10]. This 
is most easily seen by writing the complete hyperbolic cusp metric 
$g_{-1}$ on ${\mathbb R}\times T^{n-1}$ in the form
\begin{equation}\label{e2.7}
g_{-1} = r^{-2}dr^{2} + r^{2}g_{T^{n-1}}.
\end{equation}
Here $r \in (0, \infty)$ is given by $r = e^{t}$ in terms of (2.1). The 
direction $r \rightarrow 0$ gives the contracting end of the cusp, 
while the direction $r \rightarrow \infty$ gives the expanding end.

  As $r \rightarrow \infty$, the metrics $g_{BH}$ and $g_{-1}$ clearly 
approximate each other. In fact, the curvature tensor of $g_{BH}$ is 
easily calculated as follows: let $e_{i}$ be an orthonormal basis for 
$g_{BH}$ at a given point, with $e_{1}$ pointing in the $r$ direction, 
$e_{2}$ pointing in the $\theta$ direction, and $e_{i}$, $i \geq 3$ 
tangent to the toral factor. This basis diagonalizes the curvature 
tensor at every point, and the sectional curvatures $K$ in the 
corresponding 2-planes are given by
\begin{equation}\label{e2.8}
K_{12} = -1 + \frac{(n-3)(n-2)m}{r^{n-1}}, \ \ K_{1i} = -1 - 
\frac{(n-3)m}{r^{n-1}}, \ i \geq 3,
\end{equation}
$$K_{2i} = -1 - \frac{(n-3)m}{r^{n-1}}, \ i \geq 3, \ \ K_{ij} = -1 + 
\frac{2m}{r^{n-1}}, \ i, j \geq 3.$$
Thus, the curvature decays to that of the hyperbolic metric at a rate 
of $r^{-(n-1)}$, as $r \rightarrow \infty$. Let $s$ denote the geodesic 
distance to the core torus $T^{n-2}$, so that $s = s(r)$ with $ds / dr 
= V^{-1/2}$. For $r$ large, $r \sim e^{s}$, and so the curvature decays 
to $-1$ as $O(e^{-(n-1)s})$. In particular, $|W| = O(e^{-(n-1)s})$, for 
the Weyl curvature $W$. Similarly, one easily computes that 
$|\nabla^{k}R| = O(e^{-(n-1)s}) = O(r^{-(n-1)})$ for the decay of the 
covariant derivatives of the curvature tensor. 

  The function $\rho = r^{-1}$ is a smooth, geodesic defining function 
for the boundary $S^{1}\times T^{n-2} \simeq T^{n-1}$ of $D^{2}\times 
T^{n-2}$ and hence the natural conformal compactification of $g_{BH}$ 
given by 
\begin{equation}\label{e2.9}
\bar g_{BH} = \rho^{2}g_{BH},
\end{equation}
extends smoothly to the boundary to give a metric $\gamma$ on the 
conformal infinity $T^{n-1}$. Clearly, the metric $\gamma$ is the flat 
product metric $d\theta^{2} + g_{T^{n-2}}$, where the circle 
parametrized by $\theta$ has length $\beta$ given by (1.11). Note that 
the mass $m$ thus determines the length $\beta$ of the $S^{1}$ at 
conformal infinity. Further, it is important to note that (1.11) shows 
$\beta$ is strictly monotonically decreasing in $m$, $\beta'(m) < 0$. 

\bigskip

{\bf \S 2.3.} Next, we briefly discuss Einstein metrics and the 
linearization of the Einstein operator. Let $M$ be an arbitrary closed 
$n$-manifold, or the interior of a compact manifold with boundary. Let 
${\mathbb M}^{m, \alpha}$ be the space of $C^{m,\alpha}$ complete 
Riemannian metrics on $M$, i.e. complete metrics which are $C^{m, 
\alpha}$ in a smooth atlas on $M$. A more precise description of the 
topology on ${\mathbb M}^{m, \alpha}$ is given later in \S 2.4. For 
convenience, we assume $m \geq 3$, $\alpha \in (0,1)$. Similarly, let 
${\mathbb S}_{2}^{m, \alpha}$ be the space of $C^{m, \alpha}$ symmetric 
bilinear forms on $M$.

  The Einstein condition (1.1) is diffeomorphism invariant, and hence 
if $g$ is Einstein, so is $\phi^{*}g$, for any diffeomorphism $\phi$. 
In order to take this invariance into account, following Biquard [10], 
it is natural to consider the related operator
\begin{equation}\label{e2.10}
\Phi : {\mathbb M}^{m, \alpha} \longrightarrow {\mathbb S}_{2}^{m-2, 
\alpha},
\end{equation}
\begin{equation}\label{e2.11}
\Phi(g) =  Ric_{g} + (n-1)g + (\delta_{g})^*\left(\delta_{\bar g}g + 
\tfrac{1}{2} d(tr_{\bar g} g) \right).
\end{equation}
Here $\bar g$ is any fixed, (background) metric in ${\mathbb M}^{m, 
\alpha}$, $\delta$ is the divergence operator, with respect to the given 
metric, and $\delta^{*}$ is its $L^{2}$ adjoint. Recall that $\beta_{\bar g} = 
\delta_{\bar g} + \frac{1}{2}dtr_{\bar g}$ is the Bianchi operator 
associated to $\bar g$. In the applications in this paper, $\bar g$ 
will be a constructed, approximate solution to the Einstein equations, 
(called $\widetilde g$ later), while $g$ will be a metric nearby to 
$\bar g$ in the $C^{m, \alpha}$ topology. The map $\Phi$ is clearly a 
$C^{\infty}$ smooth map.

  There are two basic reasons for considering the operator $\Phi$. 
First:

\begin{lemma}\label{l2.1}
Suppose $Ric_{g} - \lambda g \leq 0$, for some $\lambda < 0$ and 
$|\beta_{\bar g}(g)|$ is bounded. If $\Phi(g) = 0$, then $g$ is 
Einstein, and 
$$Ric_{g} = -(n-1)g.$$
\end{lemma}

{\bf Proof:} This result is essentially proved in [10, Lemma I.1.4], in the 
context of asymptotically hyperbolic metrics. The proof in the case of complete 
manifolds with $Ric$ strictly negative is the same, but for completeness 
we give the proof. Applying the operator $\beta_{g}$ to both sides of 
(2.11), and using the Bianchi identity and a standard Weitzenbock formula, 
gives 
$$(D^{*}D - Ric_{g})(\beta_{\bar g}(g)) = 0.$$
Taking the inner product of this with $\beta_{\bar g}$ with respect to $g$ 
then gives $-\Delta |\beta_{\bar g}(g)|^{2} + |D\beta_{\bar g}(g)|^{2} - 
 Ric_{g}(\beta_{\bar g}(g),\beta_{\bar g}(g)) = 0$. The last two terms 
are non-negative, with the last term positive wherever $|\beta_{\bar g}(g)| > 
0$. The result then follows by a standard application of the 
maximum principle, or more precisely a maximum principle at infinity, cf. 
[35].
{\endproof}

  The map $\Phi$ is not equivariant with respect to the action of the 
diffeomorphism group, and so not every Einstein metric $h$ near $\bar g$ 
satisfies $\Phi(h) = 0$. On the other hand, the variety $\Phi^{-1}(0)$ 
gives a local slice for space of Einstein metrics near $\bar g$, transverse to 
the orbits of the diffeomorphism group, cf. [10].

  The second reason is that the form of the linearization $D\Phi$ at 
$\bar g$ has an especially simple form, cf. [10, (1.9)]:
\begin{equation}\label{e2.12}
(D_{\bar g}\Phi)(h) = \frac{1}{2}[D^{*}Dh - 2R(h) + Ric \circ h + h 
\circ Ric + 2(n-1)h].
\end{equation}
Here all metric quantities on the right are with respect to $\bar g$ and $R(h)$ 
is the action of the curvature tensor of $\bar g$ on symmetric bilinear 
forms, cf. [8, 1.131]. In particular, the operator $D_{\bar g}\Phi$ is elliptic. 
For metrics $\bar g$ of constant curvature $-1$, one easily computes that
\begin{equation}\label{e2.13}
R(h) = h - (tr h)\bar g.
\end{equation}

  For later use, we record here the Weitzenbock formula on symmetric 
bilinear forms, cf. [8, 12.69]
\begin{equation}\label{e2.14}
D^{*}Dh = (\delta d + d \delta)h + R(h) - h \circ Ric,
\end{equation}
where $d = d^{\nabla}$ is the exterior derivative induced by the metric 
connection $\nabla$, and $\delta$ is the adjoint of $d$. Hence, (2.12) 
may be rewritten in the form
\begin{equation}\label{e2.15}
2(D_{\bar g}\Phi)(h) = L(h) = (\delta d + d\delta)h - R(h) + Ric \circ 
h + 2(n-1)h.
\end{equation}
For Einstein metrics, this becomes
\begin{equation}\label{e2.16}
L(h) = (\delta d + d\delta)h - R(h) + (n-1)h.
\end{equation}
The kernel $K = {\rm Ker}L$ is the space of (essential) infinitesimal 
Einstein deformations.

\bigskip

{\bf \S 2.4.} We conclude with a discussion of topologies on the space 
of metrics that will be used below. As above ${\mathbb M}$ denotes the 
space of complete Riemannian metrics on a given manifold $M$. The 
tangent space to ${\mathbb M}$ at any point is ${\mathbb S}_{2}$ - the 
space of symmetric bilinear forms on $M$. Let ${\mathbb M}^{m}$ be the 
space of $C^m$ complete Riemannian metrics on $M$ - i.e. there exist 
(smooth) local coordinates in which the metric is $C^m$. The space 
${\mathbb M}^{m}$ may be defined intrinsically, (without use of local 
coordinates) by means of a $C^m$ norm on the tangent spaces 
$T_{g}{\mathbb M}$. Thus, given $h \in T_{g}{\mathbb M}$, define
$$||h||_{C^{m}(g)} = \sup_{x\in M}[|h|(x) + |Dh|(x) + ... + 
|D^mh|(x)],$$
where $D^{j}$ is the $j^{\rm th}$ covariant derivative; both the 
covariant derivative and (pointwise) norm are taken with respect to $g$. 
One may then define ${\mathbb M}^{m}$ to be the completion of the space of 
$C^{\infty}$ complete metrics with respect to this norm. It is standard that 
these two definitions of ${\mathbb M}^{m}$ agree.

   However, the spaces $C^m$ are not suitable for estimates for 
elliptic equations, (as in (2.12)) - which will be needed in the proof. 
For this, one must use the H\"older spaces $C^{m, \alpha}$, $\alpha \in 
(0,1)$. We are not aware of any intrinsic definition of such H\"older 
spaces of metrics and so local coordinates are needed to define them.

   For a given metric $g$ on an $n$-manifold $M$, the coordinates 
giving the optimal regularity properties for the metric are harmonic 
coordinates. Let $\rho^{m, \alpha}(x)$ be the $C^{m, \alpha}$ harmonic 
radius at $x \in M$, cf. [3]. This is the largest radius such that, for 
any $r < \rho^{m, \alpha}(x)$, the geodesic ball $B_{x}(r)$ has 
harmonic coordinates in which the metric components $g_{ij}$ satisfy
\begin{equation}\label{e2.17}
Q^{-1}(\delta_{ij}) \leq (g_{ij}) \leq Q(\delta_{ij}),
\end{equation}
\begin{equation}\label{e2.18}
\sum_{1 \leq |\beta| \leq 
m}r^{|\beta|}\sup_{y}|\partial^{\beta}g_{ij}(y)| + \sum_{|\beta| = 
m}r^{m+\alpha}\sup_{y_{1}, y_{2}} \frac{|\partial^{\beta}g_{ij}(y_{1}) 
- \partial^{\beta}g_{ij}(y_{2})|}{|y_{1} - y_{2}|^{\alpha}} \leq Q - 1.
\end{equation}
Here $Q > 1$ is a constant, fixed once and for all, (close to 1). 

 It is proved in [3] that there is a lower bound on $\rho^{m, \alpha}$, 
$\rho^{m, \alpha} \geq \rho_{0} > 0$, on any Riemannian manifold, where 
$\rho_{0}$ depends only on an upper bound for 
$||\nabla^{m-1}Ric||_{L^{\infty}}$ and a lower bound for the 
injectivity radius $inj$:
\begin{equation}\label{e2.19}
||\nabla^{m-1}Ric||_{L^{\infty}} \leq \Lambda < \infty, \ {\rm and} \ 
inj \geq i_{0} > 0.
\end{equation} 

  Given a Riemannian manifold $(M, g)$ satisfying (2.19), choose a 
covering ${\mathcal U}_{\lambda}$ of $(M, g)$ by a collection of 
$\rho_{0}/2$ balls such that the $\rho_{0}/4$ balls are disjoint. The 
bounds (2.17) imply a uniform upper bound on the multiplicity of such a 
covering. Now let $g'$ be another metric on $M$ and set $g' - g = h$, 
so that $h \in {\mathbb S}_{2}(M)$. As in (2.18), define then
\begin{equation}\label{2.20}
||g'||_{C^{m, \alpha}(g)} \equiv ||h||_{C^{m, \alpha}} = \sup_{\lambda} 
\{ \sum_{1 \leq |\beta| \leq 
m}\rho_{0}^{|\beta|}\sup_{y}|\partial^{\beta}h_{ij}^{\lambda}(y)| + 
\sum_{|\beta| = m}\rho_{0}^{m+\alpha}\sup_{y_{1},y_{2}} 
\frac{|\partial^{\beta}h_{ij}^{\lambda}(y_{1}) - 
\partial^{\beta}h_{ij}^{\lambda}(y_{2})|}{|y_{1} - y_{2}|^{\alpha}} \},
\end{equation}
where the components $h_{ij}^{\lambda}$ are taken in local $g$-harmonic 
coordinates $u_{i}^{\lambda}$ satisfying (2.17)-(2.18), and the 
supremum (2.20) is taken over all such local coordinate systems in 
${\mathcal U}_{\lambda}$. 
  
 This defines the $C^{m, \alpha}$ topology on ${\mathbb M}$, denoted as 
${\mathbb M}^{m, \alpha}$, in a neighborhood of a given metric $g$ on 
which one has bounds on the Ricci curvature and injectivity radius as 
above. 

\medskip

  In the course of the arguments to follow, we will have good control 
on the Ricci curvature, to all orders. However, for the classes of 
metrics to be considered, there will not be a uniform lower bound on 
the injectivity radius; this will cause the norm (2.20) to degenerate.

   In general, when the injectivity radius is very small, the geometry 
of small balls may be very complicated; (this involves the structure of 
collapsed manifolds in the sense of Cheeger-Gromov, with bounds on 
Ricci curvature). Fortunately, we need only deal with situations where 
the metrics have {\it bounded local covering geometry}, in the 
following sense.

\begin{definition}\label{d2.2}
Let $i_{0} > 0$ be given. Then $(M, g)$ has bounded local covering 
geometry, (with respect to $i_{0}$), if for any $x$ where $inj(x) \leq i_{0}$, 
there is a finite covering space $\bar B_{x, i_{0}}$ of the geodesic 
ball $B_{x}(i_{0})$ with $diam_{g} \bar B_{x, i_{0}} \leq 1$ and 
$$inj_{g}(\bar x) \geq i_{0}.$$
Here $\bar x$ is a lift of $x$ to $\bar B_{x, i_{0}}$, and $g$ is 
lifted to $\bar B_{x, i_{0}}$ so that the projection is a local 
isometry.
\end{definition} 

Thus, by passing to a finite covering space locally, one can unwrap to 
obtain a metric of bounded geometry, and thus good local harmonic 
coordinates as in (2.17)-(2.18), given suitable control on the Ricci 
curvature. The degree of the covering of course depends on the 
injectivity radius at $x$; the smaller the injectivity radius, the 
larger the degree of the covering. This definition depends on a choice 
of $i_{0}$. For our purposes, $i_{0}$ will be a fixed small number, 
depending only on dimension, throughout the paper. One may take for 
instance $i_{0}$ to be a fixed small multiple of the Margulis constant 
in dimension $n$, cf. [18], [31].

  Let $(M, g)$ be any complete Riemannian manifold satisfying the bound 
\begin{equation}\label{e2.21}
||\nabla^{m-1}Ric||_{L^{\infty}} \leq \Lambda < \infty,
\end{equation}
and which has bounded local covering geometry with respect to $i_{0}$. 
One may then define a ``modified'' $C^{m, \alpha}$ norm  $\widetilde C^{m, 
\alpha}$ of a metric $g'$ by setting $h = g' - g$, and defining 
\begin{equation}\label{e2.22}
||g'||_{\widetilde C^{m, \alpha}(g)} \equiv ||h||_{\widetilde C^{m, 
\alpha}}
\end{equation}
exactly as in (2.20) where the charts are defined in finite covering 
spaces as above in regions where the injectivity radius is $\leq i_{0}$.

\section{Proof of Theorem 1.1.}
\setcounter{equation}{0}

  This section is mainly concerned with the proof of Theorem 1.1. 
Following the proof, Proposition 3.9 proves that the homeomorphism type of 
$M_{\bar \sigma}$ is determined up to finite ambiguity by the curves in 
$\bar \sigma$. Corollary 3.11 is a version of Theorem 1.1 on non-toral 
ends.

\medskip

  We break the proof of Theorem 1.1 into two main steps.

\medskip

{\bf Step I. Construction of the Approximate Solution.}

\medskip

  One begins with a complete non-compact hyperbolic $n$-manifold $(N, 
g_{-1})$ of finite volume, and its collection of toral ends 
$T^{n-1} \times {\mathbb R}^{+}$. Fix any such end $E$, and a flat 
torus $T^{n-1} \subset E$, normalized as in \S 2.1. Given a simple 
closed geodesic $\sigma$ in $(T^{n-1}, g_{0})$, the discussion in \S 
2.1 describes Dehn filling topologically on the end $E$. In this step, 
we construct this filling metrically.

   Consider the standard toral AdS black hole metric (2.6) on 
$D^{2}\times T^{n-2}$:
\begin{equation}\label{e3.1}
g_{BH} = V^{-1}dr^{2} + Vd\theta^{2} + r^{2}g_{T^{n-2}}.
\end{equation}
On the universal cover $D^{2}\times {\mathbb R}^{n-2}$, the metric 
$g_{BH}$ lifts to a metric $\widetilde g_{BH}$ of the form (3.1), with 
flat metric on $T^{n-2}$ lifted to ${\mathbb R}^{n-2}$, i.e.
\begin{equation}\label{e3.2}
\widetilde g_{BH} = V^{-1}dr^{2} + Vd\theta^{2} + r^{2}(ds_{1}^{2} + 
... +ds_{n-2}^{2}).
\end{equation}
The change of variable $r \rightarrow r_{m} = m^{1/(n-3)}r$ shows that 
the metrics $\widetilde g_{BH} = \widetilde g_{BH}(m)$ are all 
isometric. Thus, for convenience, we fix $m$ once and for all, by 
setting, (for example), $m = \frac{1}{2}$, so that $r_{+} = 1$. 

  Let $D(R) = \{r \leq R\}$ in $(D^{2}\times {\mathbb R}^{n-2}, 
\widetilde g_{BH})$ and let $S(R) = \partial D(R) = \{r = R\}$. The 
induced metric on the boundary $S(R)$ is then a flat metric
\begin{equation}\label{e3.3}
V(R)d\theta^{2} + (dt_{1}^{2} + ... +dt_{n-2}^{2})
\end{equation}
on $S^{1}\times {\mathbb R}^{n-2}$, where $t_{i} = R_{i}s_{i}$ are 
coordinates on ${\mathbb R}^{n-2}$. Choose $R$ so that 
\begin{equation}\label{e3.4}
V(R)^{1/2} \cdot \beta = L(\sigma).
\end{equation}
Thus, the length of $S^{1} \times \{pt\} \subset S(R)$ equals 
$L(\sigma)$. Recall that $V = V_{m}$ and $\beta = \beta(m)$ are 
determined since $m = \frac{1}{2}$. 

  Given the flat structure $g_{0}$ on the torus $T^{n-1}$, observe that 
there is a unique (up to conjugacy) free isometric ${\mathbb Z}^{n-2}$ 
action on the flat product $\partial S(R) = S^{1}\times {\mathbb 
R}^{n-2}$ such that the projection map to the orbit space
\begin{equation}\label{e3.5}
\pi:S^{1}\times {\mathbb R}^{n-2} \rightarrow T^{n-1}
\end{equation}
satisfies $\pi(S^{1}) = \sigma$, and the flat structure on $T^{n-1}$ 
induced by $\pi$ is the given $g_{0}$. In fact the map $\pi$ is just 
the covering space of $(T^{n-1}, g_{0})$ corresponding to the subgroup  
$\langle \sigma \rangle \subset \pi_{1}(T^{n-1})$. In more detail, 
$\sigma = \sum \sigma^{i} v_{i}$ may be viewed as a vector in ${\mathbb 
R}^{n-1}$. This may be completed to an integral  basis $(\sigma, b_{2}, 
..., b_{n-1})$ of ${\mathbb R}^{n-1}$ in such a way that the lattice 
generated by $(\sigma, b_{2}, ..., b_{n-1})$ equals the lattice 
generated by $(v_{1}, ..., v_{n-1})$, i.e. there is a matrix in 
$SL(n-1, {\mathbb Z})$ taking $(v_{1}, ..., v_{n-1})$ to $(\sigma, 
b_{2}, ..., b_{n-1})$. Without loss of generality, we may assume that 
the length of the projection of each $b_{i}$ onto $\sigma$ has length 
at most $|\sigma|$, i.e. $|\langle b_{i}, \sigma \rangle | < 
|\sigma|^{2}$. Then $S(R)$ may be identified with ${\mathbb R}^{n-1} / 
\langle \sigma \rangle$, where $\langle \sigma \rangle \simeq {\mathbb 
Z}$ is the group generated by $\sigma$. The vectors $b_{2}, ..., 
b_{n-1}$ generate a ${\mathbb Z}^{n-2}$ action on ${\mathbb R}^{n-1}$ 
commuting with $\langle \sigma \rangle$, and hence generate a ${\mathbb 
Z}^{n-2}$ action on $S(R)$. The map $\pi$ is then the map to the orbit 
space of this action.

  This ${\mathbb Z}^{n-2}$ action extends radially to an isometric 
action on the domain $D(R)$ contained in the universal cover 
$D^{2}\times {\mathbb R}^{n-2}$. To see this, the isometry group of 
$\widetilde g_{BH}$ is Isom$(S^{1})\times$ Isom$({\mathbb R}^{n-2})$, 
corresponding to rotations in the $\theta$-circle and Euclidean 
isometries on ${\mathbb R}^{n-2}$. Any isometry of the boundary 
$\partial D(R) = S(R)$ thus extends uniquely to an isometry of $D(R)$. 
It is clear that the resulting action on $D(R)$ or the full universal 
cover $D^{2} \times {\mathbb R}^{n-2}$ is smooth and free.

  The quotient space $(D^{2} \times {\mathbb R}^{n-2}) / {\mathbb 
Z}^{n-2} \simeq D^{2} \times T^{n-2}$ gives the (twisted) toral AdS 
black hole metric 
\begin{equation}\label{e3.6}
g_{BH} = [V^{-1}dr^{2} + Vd\theta^{2} + r^{2}g_{{\mathbb R}^{n-2}}] / 
{\mathbb Z}^{n-2}.
\end{equation}
If now $D(R)$ denotes the domain $\{r \leq R\}$ in the quotient space, 
the boundary $S(R) = \partial D(R)$ is isometric to the initially given 
flat torus $(T^{n-1}, g_{0})$. 

  As $r$ varies over $(r^{+}, R]$, the tori $S(r)$ with metric induced 
from $g_{BH}$ give a curve of flat metrics on $T^{n-1}$. To describe 
this curve, let $\lambda(r) = \beta V^{1/2}(r)/|\sigma| = 
(V(r)/V(R))^{1/2}$, so that $\lambda(r) \in (0,1]$. Then the torus 
$S(r)$ is generated by $(\sigma(r), b_{2}(r), ..., b_{n-1}(r))$, where 
\begin{equation}\label{e3.7}
\sigma(r) = \lambda(r)\cdot \sigma, \ {\rm and} \ b_{i}(r) = b_{i} + 
(\lambda(r) - 1)(\langle b_{i}, \sigma \rangle /|\sigma|^{2})\sigma.
\end{equation}
Note that $L(\sigma(r)) \rightarrow 0$, as $r \rightarrow r_{+}$, and 
at $\{r = r^{+}\}$, the generators $b_{i}(r^{+})$ of the core 
$(n-2)$-torus $T^{n-2}$ are orthogonal to $\sigma$. 

  Observe also that for $R$ large, equivalently $L(\sigma)$ large, the 
core totally geodesic $T^{n-2}$ at $r = r_{+}$ shrinks to 0 size; in 
fact 
$$diam T^{n-2} \sim R^{-1}.$$
In particular, the injectivity radius of $g_{BH}$ at and near $T^{n-2}$ 
is $O(R^{-1})$. On the other hand, the metrics $g_{BH}$ clearly have 
uniformly locally bounded covering geometry, independent of $R$, cf. \S 
2.4. When $n = 3$, the metric $g_{BH}$ is hyperbolic, and is a complete 
hyperbolic tube metric about a closed geodesic of length $\sim R^{-1}$, 
cf. [18]. 

\medskip

  Since the boundaries $(\partial D(R), g_{BH}) = S(R)$ and $(T^{n-1}, 
g_{0}) \subset (E, g_{-1})$ are isometric, they may be identified; this 
gives the Dehn filling $M_{\sigma}$ of the end $E$ along the curve 
$\sigma$.

  Although the intrinsic flat metrics on the boundaries agree, the union 
of the two ambient metrics $g_{BH}$ and $g_{-1}$ forms a corner at the seam 
$\partial S(R)$. To estimate the difference of the metrics, it is convenient 
to write the hyperbolic cusp metric $g_{-1}$ from (2.7) in the form
\begin{equation}\label{e3.8}
g_{-1} = r^{-2}dr^{2} + r^{2}g_{\frac{1}{R}T^{n-1}},
\end{equation}
so that $R^{2}g_{R^{-1}T^{n-1}} = g_{0}$. This just amounts to 
replacing $r$ by $r/R$ in (2.7) and has the effect that the glueing 
seam is located at $\{r = R\}$ for both metrics. Thus, comparing (3.6) 
and (3.8), one sees that $g_{BH}$ and $g_{-1}$ differ on the order of 
$O(R^{1-n})$ near the seam. A simple computation also shows that the 
$2^{\rm nd}$ fundamental forms $A_{-1}$ and $A_{BH}$ of the boundary 
with respect to $g_{-1}$ and $g_{BH}$ are 
$$A_{-1} = g_{-1}|_{T^{n-1}},$$
$$A_{BH} \sim (1 + O(R^{1-n}))g_{BH}|_{T^{n-1}}.$$
Thus, the $2^{\rm nd}$ fundamental forms differ on the order of 
$O(R^{1-n})$. Similarly, from (2.8), the curvatures of the two metrics 
also differ on the order of $O(R^{1-n})$. 

  One may then smooth the corner at the toral seam $S(R)$ by setting
\begin{equation}\label{e3.9}
\widetilde g = [\widetilde V^{-1}dr^{2} + \widetilde V d\theta^{2} + 
r^{2}g_{{\mathbb R}^{n-2}}] / {\mathbb Z}^{n-2},
\end{equation}
where, recalling $m = \frac{1}{2}$,
$$\widetilde V = r^{2} - \frac{\chi \circ r}{r^{n-3}}.$$
Here $\chi: {\mathbb R} \rightarrow {\mathbb R}$ is a smooth function 
satisfying $\chi(r) = 1$, for $1 \leq r \leq R/2$, $\chi(r) = 0$, for 
$r \geq 2R$ and $|\partial^{k}\chi| = O(R^{-k})$. Note here also that 
the geodesic distance between the $r$-levels $R/2$ and $2R$ is on the 
order of $1$. 

  The smooth metric $\widetilde g$ extends to a globally defined metric 
on $M_{\sigma}$, by letting $\widetilde g$ be the hyperbolic metric on 
$N$. This process may be carried out on any collection of toral ends 
$E_{j}$, $1 \leq j \leq p$ of $N$ and gives a smooth metric $\widetilde 
g$ on $M_{\bar \sigma} = M(\sigma_{1}, ..., \sigma_{p})$. This gives a 
collection of numbers $\bar R = (R_{1}, ..., R_{p})$ corresponding to 
$\{\sigma_{j}\}$ via (3.4). Let 
$$R_{min} = \min_{j}R_{j} .$$
We also set $R_{max} = \max_{j}R_{j}$, but note that $R_{max} = \infty$ if 
$p < q$, i.e. if there is an end of $N$ which is not capped off by Dehn filling. 

These metrics will be called approximate solutions of the Einstein equation (1.1). 

   The discussion above proves the following result:

\begin{proposition}\label{p3.1}
The approximate solutions $\widetilde g$ constructed above on $M_{\bar 
\sigma}$ are complete, and of uniformly bounded local covering 
geometry. Outside a tubular neighborhood $U_{j}$ of radius 1 
about each fixed torus $T_{j}^{n-1}$, $1 \leq j \leq p$, $\widetilde g$ is 
the hyperbolic metric $g_{-1}$ on $N$ or the black hole metric (3.6) on 
$D^{2}\times T^{n-2}$. The curvature of $\widetilde g$ is uniformly bounded 
by that of $g_{BH}$, in that its sectional curvature is bounded by the values 
in (2.8) with $r = 2m = 1$; if $n = 3$, then the curvature of $\widetilde g$ 
is $-1 + O(R_{min}^{-2})$. 

  The metric $\widetilde g$ satisfies the Einstein equation
\begin{equation}\label{e3.10}
Ric_{\widetilde g} + (n-1)\widetilde g = 0,
\end{equation}
outside $U = \cup U_{j}$, while inside each $U_{j}$,
\begin{equation}\label{e3.11}
Ric_{\widetilde g} + (n-1)\widetilde g = O(R_{j}^{1-n}), \ {\rm and} \ 
|\nabla^{k}Ric_{\widetilde g}| = O(R_{j}^{1-n}), 
\ {\rm for\  any} \ k < \infty.
\end{equation}
\end{proposition}
{\endproof}

\smallskip

{\bf Step II. Analysis of the Linearization.}
 
\medskip

 The strategy now is to use the inverse function theorem to perturb the 
approximate solution $\widetilde g$ constructed on $M = M_{\bar \sigma}$ 
into an exact solution of the Einstein equation (1.1). To do 
this, one needs to study the linearization of the Einstein operator 
(2.11) at $\widetilde g$. Thus, set 
$$L = 2D_{\widetilde g}\Phi,$$
so that, from (2.12),
\begin{equation}\label{e3.12}
L(h) = D^{*}Dh - 2R(h) + Ric\circ h + h\circ Ric + 2(n-1)h. 
\end{equation}
where the metric quantities on the right are with respect to $\widetilde g$. 
For reasons that will soon be apparent, we assume throughout Step II that 
\begin{equation} \label{e3.13}
M = M_{\bar \sigma} \ {\rm is} \ {\rm compact ,}
\end{equation}
so that $p = q$ and all ends of $N$ are Dehn filled. This assumption will 
be removed later, cf.  . Under the assumption (3.13), we will show that 
$L$ is invertible on suitable function spaces, and obtain a bound on the 
inverse $L^{-1}$, for all sufficiently large Dehn fillings $\bar \sigma$. 
In addition, these statements hold for metrics sufficiently close to 
$\widetilde g$.

\medskip

 To begin, as function spaces, we will use the modified H\"older spaces 
and norms, discussed in \S 2.4; these are well-adapted to the 
approximate solutions $\widetilde g$, since by (3.11), the metrics 
$\widetilde g$ have uniformly bounded Ricci curvature, (in fact 
uniformly bounded curvature), to all orders, for all $\bar \sigma$. 
Further, the metrics $\widetilde g$ have uniformly bounded local 
covering geometry, again independent of $\bar \sigma$.

  Thus, fix any $m \geq 3$, $\alpha\in (0,1)$. The map $\Phi$ is a smooth map
$$\Phi : {\mathbb M}^{m,\alpha} \rightarrow  {\mathbb S}_{2}^{m-2,\alpha},$$ 
with derivative at $\widetilde g$, (modulo the factor of 2), a 
smooth linear map
$$L: {\mathbb S}_{2}^{m,\alpha} \rightarrow  {\mathbb S}_{2}^{m-2,\alpha},$$
\begin{equation}\label{e3.14}
L(h) = f.
\end{equation}
Recall that $R_{max} = \max_{j}R_{j}$. 
\begin{proposition}\label{p3.2}
 For $M = M_{\bar \sigma}$ as in (3.13) with $\bar \sigma$ sufficiently large, 
there is a constant $\Lambda$, independent of $\bar \sigma$, such that
\begin{equation}\label{e3.15}
||h||_{\widetilde C^{m,\alpha}} \leq  \Lambda (\log R_{max})
||L(h)||_{\widetilde C^{m-2,\alpha}}. 
\end{equation}
It follows that $L$ is invertible and the norm of $L^{-1}: 
{\mathbb S}_{2}^{m-2,\alpha} \rightarrow  {\mathbb S}_{2}^{m,\alpha}$ is 
uniformly bounded by $\Lambda \log R_{max}$.
\end{proposition}

{\bf Proof:}
 Note first that the estimate (3.15) is local, in the sense that the norms are 
taken with respect to controlled local harmonic coordinate charts (2.17)-(2.18), 
in suitable covers where the injectivity radius is small. 

 The operator $L$ is an elliptic operator on $h$, and by an examination 
of the form of $L$ in (3.12), one has uniform control on all the 
coefficients of $L$ in local harmonic coordinates. More precisely, the 
leading order term $D^{*}D$ has (uniformly bounded) $C^{m, \alpha}$ 
coefficients, while the $0$-order terms involving curvature give 
(uniformly bounded) $C^{m-2, \alpha}$ coefficients. Hence, the Schauder 
estimates for elliptic systems, cf. [17], [27], give the estimate
\begin{equation}\label{e3.16}
||h||_{\widetilde C^{m,\alpha}} \leq  \Lambda\{||L(h)||_{\widetilde 
C^{m-2,\alpha}}+ ||h||_{L^{\infty}}\}, 
\end{equation}
where $\Lambda$ is independent of the Dehn filling. Note that the 
$L^{\infty}$ norm is invariant under passing to (local) covering 
spaces. Setting $f = L(h)$ as above, it then suffices to prove that 
there exists $\Lambda < \infty$ such that
\begin{equation}\label{e3.17}
||h||_{L^{\infty}} \leq  \Lambda \log R_{max}||f||_{\widetilde C^{m-2,\alpha}}. 
\end{equation}

 The claim is that the estimate (3.17) holds provided all Dehn fillings 
$\sigma_{j} \in \bar \sigma$ are sufficiently large with $\Lambda$ independent 
of $\bar \sigma$. We prove this by contradiction; some comments on the 
possibility of a more effective proof are given in Remark 3.5 below.

 Thus, suppose (3.17) is false. Then there is a sequence of Dehn-filled 
manifolds $M_{i} = M_{\bar \sigma_{i}}$, with $(\sigma_{j})_{i} 
\rightarrow  \infty$ for each $(\sigma_{j})_{i} \in \bar \sigma_{i}$, 
together with approximate solutions $\widetilde g_{i}$ on $M_{i}$, and 
symmetric forms $h_{i}\in {\mathbb S}_{2}^{m,\alpha}(M_{i})$, such that
\begin{equation}\label{e3.18}
||h_{i}||_{L^{\infty}} = 1, \ \ {\rm but} \ \ \log(R_{max})_{i} 
||f_{i}||_{\widetilde C^{m-2,\alpha}} \rightarrow  0, 
\end{equation}
where $f_{i} = L_{i}(h_{i})$. Observe that the estimate (3.16) now 
implies that 
\begin{equation}\label{e3.19}
||h_{i}||_{\widetilde C^{m, \alpha}} \leq \Lambda,
\end{equation}
where $\Lambda$ is fixed, (independent of $i$).

  The idea of the proof then is to pass to limits, and produce a 
non-trivial limit form $h$ in Ker $L$. Roughly speaking, the 
manifold $(M_{i}, \widetilde g_{i})$ divides into three regions - the 
hyperbolic region $N$, the cusp regions and the black hole regions. The cusp 
regions arise as a transition between the hyperbolic and black hole geometries. 
A well-known argument, essentially due to Calabi [14], implies that $L$ 
has no kernel on $N$. We will prove that the cusp and black hole regions 
also have no kernel. Taken together, these facts will give a contradiction 
to the behavior (3.18). We now supply the details of this description.

\medskip

  First, we prove an elementary Lemma, (which will be needed however only in 
Appendix A). 
\begin{lemma}\label{l3.3}
Under the assumptions (3.18), one has
\begin{equation}\label{e3.20}
||tr h_{i}||_{L^{\infty}} \rightarrow 0 \ \ {\rm as} \ i \rightarrow \infty.
\end{equation}
\end{lemma}
{\bf Proof:} Taking the trace of (3.16), using (3.12) and the fact that 
$tr R(h) = \langle Ric, h \rangle$, gives, (dropping the $i$ from the notation),
$$-\Delta tr h - \frac{2s}{n}tr h = tr f + \langle z, h \rangle,$$
where $z$ is the trace-free Ricci curvature. The metric $\widetilde g$ is 
almost Einstein; $|z| \leq O([R_{min}]_{i}^{-(n-1)})$, cf. (3.11). 
Since $|h|$ is uniformly bounded, one has $|\langle z, h \rangle| \rightarrow 0$, 
as $i \rightarrow \infty$. Since also $|f| \rightarrow 0$, the right side of the 
equation above tends to 0 in $L^{\infty}$ as $i \rightarrow \infty$. The left side 
is a strictly positive operator, since $s \sim -n(n-1)$. Hence, the result follows 
by evaluating the equation above at points realizing the maximum and minimum of 
$tr h$. 
{\endproof}

  We now continue with the proof of Proposition 3.2 itself. Let $T = 
\cup T_{j}^{n-1}$ be the collection of tori $T^{n-1}$ in $N$ to which 
the solid tori are attached by Dehn filling, and let $N_{T}$ be the 
hyperbolic manifold obtained by removing these cusp ends 
$T_{j}^{n-1}\times {\mathbb R}^{+}$ from $N$. The manifold $M_{i} = 
M_{{\bar \sigma}_{i}}$ is a union of black hole and hyperbolic regions:
$$M_{i} = \{\cup_{j} D(R_{i}^{j})\} \cup N_{T},$$
where for each $j$, $D(R_{i}^{j})$ is the black hole region defined as 
following (3.6); thus $\partial D(R_{i}^{j})$ is attached to 
$T_{j}^{n-1}$. Observe that for any fixed $j$, $R_{i}^{j} \rightarrow \infty$, 
as $i \rightarrow \infty$. In the following, we will often work with each component 
of $D(R_{i}^{j})$ separately, and thus usually drop $j$ from the notation. 

  Let $x_{i}$ be a sequence of base points in $(M_{i}, \widetilde g_{i})$. By 
passing to a subsequence if necessary, we may assume that $\{x_{i}\}$ has 
exactly one of the following behaviors:

  (i). (Hyperbolic) One has 
\begin{equation} \label{e3.21}
dist_{\widetilde g_{i}}(x_{i}, y_{0}) < \infty, 
\end{equation}
for some fixed point $y_{0} \in N$. In this case, the pointed sequence 
$(M_{i}, \widetilde g_{i}, x_{i})$ converges in the pointed Gromov-Hausdorff 
topology, and smoothly and uniformly on compact sets, to the limit $(N, g_{-1}, x)$, 
$x = \lim x_{i}$; $(N, g_{-1})$ is the original hyperbolic manifold.  

  (ii). (Cusps) For all $j$, 
\begin{equation} \label{e3.22}
dist_{\widetilde g_{i}}(x_{i}, (T_{j}^{n-2})_{i}) \rightarrow 
\infty, \ {\rm and} \  dist_{\widetilde g_{i}}(x_{i}, y_{0}) \rightarrow \infty,  
\end{equation}
where $T_{j}^{n-2}$ is the core torus of the Dehn filling on $E_{j}$, $1 \leq j \leq q$. 
In this case, the pointed sequence $(M_{i}, \widetilde g_{i}, x_{i})$ collapses. 
However, as discussed below, one may unwrap the collapse and obtain a complete 
limit which is a complete hyperbolic cusp as in (2.7). 

  (iii). (Black hole) For some $j$, 
\begin{equation} \label{e3.23}
dist_{\widetilde g_{i}}(x_{i}, (T_{j}^{n-2})_{i}) < \infty.
\end{equation}
 Again the pointed sequence $(M_{i}, \widetilde g_{i}, x_{i})$ collapses, but 
by passing to a subsequence, the collapse may be unwrapped and one obtains 
convergence to a complete black hole metric (2.6). 

\medskip

   We deal with each of these cases in turn. 

  {\bf Case (i)}. The forms $h_{i}$ satisfying (3.18) converge smoothly (in a 
subsequence) to a limit form $h$ on the complete manifold $N$ satisfying 
\begin{equation} \label{e3.24}
L(h) = 0,
\end{equation}
i.e. $h$ is an infinitesimal Einstein deformation of the hyperbolic metric on $N$. 

   Now we use the form (2.15) for the linearization $L = 2D\Phi$: recall 
this is
$$L(h) = (\delta d + d\delta)h - R(h) + Ric \circ h + 2(n-1)h .$$
Since $N$ is hyperbolic, $R(h) = h - (trh)g$, and so 
$$L(h) = (\delta d + d\delta)h - R(h) + (n-1)h .$$

  Pick any $r_{0}$ large and pair (3.24) with $h$. Integrating by parts 
over the domain $N_{r_{0}} = \{r \geq r_{0} > 0\}$, where $r$ is the parameter 
for any of the cusp ends of $N$ as in (2.7), one thus obtains
$$\int_{N_{r_{0}}}|dh|^{2} + |\delta h|^{2} + (n-2)|h|^{2} + (tr h)^{2} = 
\int_{\partial N_{r_{0}}}Q(h, \partial h),$$
where the boundary term involves only $h$ and its first derivative. By 
(3.18) and (3.19), $Q$ is thus uniformly bounded, while the volume form of 
$\partial N_{r_{0}}$ is $O(e^{-(n-1)r_{0}})$. Letting $r_{0} \rightarrow 0$, 
it follows that
$$h \equiv 0 \ {\rm on} \ N.$$
By the smooth convergence of $h_{i}$ to the limit form $h$, it follows that 
$h_{i}(x_{i}) \rightarrow 0$ for $x_{i}$ satisfying (3.21). This shows that 
in fact there is an exhaustion $K_{j} \subset N$, with $K_{j} \subset M_{i}$ 
for $i = i(j)$ sufficiently large, a sequence $\varepsilon_{j} \rightarrow 0$, 
and a subsequence $\{h_{i_{j}}\}$ of $\{h_{i}\}$ such that 
\begin{equation} \label{e3.25}
|h_{i_{j}}(x)| \leq \varepsilon_{j} \ \ \forall x \in K_{j} .
\end{equation}
In the following, we work only with this subsequence, and relabel 
$\{h_{i_{j}}\}$ to $\{h_{i}\}$. This shows that the support of $h_{i}$ must 
either wander down the cusp-like regions of $(M_{i}, \widetilde g_{i})$, or 
meet the black hole region of $(M_{i}, \widetilde g_{i})$. 

\medskip

  {\bf Case (ii).} 

  In this case, $x_{i}$ becomes further and further distant from any given point in 
$N$, as well as any of the black hole regions. Without loss of generality, assume that 
$\{x_{i}\}$ is contained in a fixed end $E$ of $N$. Then (3.22) is equivalent to the 
statements that $(r/R_{i})(x_{i}) \rightarrow 0$, and $r(x_{i}) \rightarrow \infty$ 
as $i \rightarrow \infty$. 
 
  By construction, the manifolds $(M_{i}, \widetilde g_{i}, x_{i})$ are collapsing 
in domains of uniformly bounded diameter about $x_{i}$. However, this collapse may 
be unwrapped, (cf. \S 2.4 and Proposition 3.1), in larger and larger finite covering 
spaces to obtain a complete limit manifold $(C, g_{-1}, x)$. The limit is clearly the 
complete hyperbolic cusp metric (2.7) on ${\mathbb R}\times T^{n-1}$, with parameter 
$r$ normalized so that $r(x) = 1$. Similarly, the forms $h_{i}$, when lifted to forms 
$\widetilde h_{i}$ on the covering spaces, are uniformly bounded in 
$\widetilde C^{m, \alpha}$. Hence, a subsequence converges in the 
$\widetilde C^{m, \alpha'}$ topology, for any $\alpha' < \alpha$, to a limit 
form $\widetilde h$ satisfying, by (3.18),
\begin{equation}\label{e3.26}
L(\widetilde h) = D^{*}D\widetilde h - 2R(\widetilde h) = 0,
\end{equation}
on $(C, g_{-1})$, i.e. $\widetilde h$ is an infinitesimal Einstein deformation. 
Since the forms $\widetilde h_{i}$ have been lifted to covering spaces, they are 
invariant under the corresponding group of covering transformations. These groups 
restrict to cyclic groups ${\mathbb Z}_{k_{i}}$ acting on each circle $S^{1}$ in 
$T^{n-1} = S^{1}\times S^{1} \cdots \times S^{1}$, with $k_{i} \rightarrow \infty$ 
as $i \rightarrow \infty$. As $i \rightarrow \infty$, these covering groups converge 
to the isometric $T^{n-1}$ action on $(C, g_{-1})$. Hence, by the smooth convergence, 
the limit form $\widetilde h$ is also $T^{n-1}$ invariant. This implies that 
$\widetilde h$ has the form 
\begin{equation} \label{e3.27}
\widetilde h = \sum h_{ab}(r)\theta^{a} \cdot \theta^{b},
\end{equation}
where $h_{ab}$ is a function of $r$ only, and $\theta^{a}$ is the natural 
orthonormal coframing of the cusp metric (2.7), with $\theta^{1} = r^{-1}dr$. 
It is also clear that $\widetilde h$ is bounded on the complete cusp $C$, since the 
bound (3.18) on $h$ passes continuously to the limit by (3.19). 

  It is shown in Appendix A  that $\widetilde h$ then necessarily satisfies 
$h_{1a} = 0$ for any $a$, (see (A.11) and (A.13)), while for any $a,b \geq 2$, 
the coefficient functions $h_{ab}$ satisfy
\begin{equation} \label{e3.28}
\Delta h_{ab} = r^{2}h_{ab}'' + nrh_{ab} = 0,
\end{equation}
see (A.8). Here $r \in (0, \infty)$ and again $r(x) = 1$. (The proof of these 
statements is deferred to Appendix A, since it is purely computational, and 
unrelated to the issues at hand). The general solution of (3.28) is given 
by $c_{1}r^{-(n-1)} + c_{2}$, cf. (A.9). Since $\widetilde h$ is bounded 
on $C$, it follows that
\begin{equation} \label{e3.29}
h_{ab} = const = c_{ab}.
\end{equation}
Geometrically, this means that all bounded $T^{n-1}$-invariant infinitesimal 
Einstein deformations of the cusp metric arise from deformations of the flat 
structure on $T^{n-1}$. 

  However, the constants $c_{ab}$ in (3.29) may apriori vary with different 
choices of the base point sequence $\{x_{i}\}$. (For instance, consider the 
function $q(r) = \sin(\log r)$; any sequence $r_{i} \rightarrow \infty$ has 
a subsequence such that $q(r)$ converges to a constant on 
$[-k+r_{i}, k+r_{i}]$, for any given $k$. Nevertheless, the constants vary 
with different choices of sequence $r_{i}$). 

  We claim that all constants $c_{ab}$ in (3.29) satisfy
\begin{equation} \label{e3.30}
c_{ab} = 0,
\end{equation}
for all $x_{i}$ satisfying (3.22). The proof of (3.30) requires the 
assumption (3.15), not just the weaker the assumption that 
$||f_{i}||_{\widetilde C^{m-2,\alpha}} \rightarrow 0$. 

  To prove (3.30), return to the black hole metric (3.6), viewed as 
part of the approximate solution $\widetilde g = \widetilde g_{i}$. The 
injectivity radius and diameter of the tori $T^{n-1}(r)$ then satisfy 
$inj(T^{n-1}(r)) \sim O(r/R)$ and $diam(T^{n-1}(r)) \sim O(r/R)$; recall here 
that $R = R_{i}^{j} \rightarrow \infty$, as $i \rightarrow \infty$, for any 
given $j$. To see this, as discussed in \S 2.2, the parameter $r$ and the 
geodesic distance $s$ from the black hole horizon are related by $r \sim e^{s}$, 
for $r$ large. Let $R = e^{S}$. Then the diameter and injectivity radius of 
the torus at the locus $r$ are approximately $e^{s-S} \sim r/R$, as claimed. 

   As above, we then unwrap in large covering spaces so that $inj(T^{n-2}) \sim 1$, 
and $diam(T^{n-2}) \sim 1$. The lifted forms $h = h_{i}$ are then invariant under 
the corresponding covering transformations; here and in the following, we drop the 
tilde from the notation. Given any fixed, large $i$ and with $h = h_{i}$, let 
$$h_{ab}(r) = \frac{1}{vol T^{n-1}(r)}\int_{T^{n-1}(r)}h_{ab}(r,\theta)d\theta$$
be the average of $h_{ab}$ over $T^{n-1}(r)$. The same definition applies to 
$f_{ab}(r)$, so that $h(r)$, $f(r)$ are $T^{n-1}$-invariant forms, as in (3.27). 
Abusing notation slighly, let $U(r) = \{x \in E: r(x) \in [\frac{1}{2}r, 2r]\}$, 
so that $U(r)$ is a tubular neighborhood about $T^{n-1}(r)$ of geodesic size on 
the order of 1, independent of $r$. Using (3.18)-(3.19), we note that one has
$$||h - h(r)||_{C^{2}(U(r))} = O(\frac{r}{R}) \ {\rm and} \ 
||f - f(r)||_{C^{0}(U(r))} = O(\frac{r}{R}),$$
independent of $i$. This is because the coefficients of the lifted forms 
$h = h_{i}$ and $f = f_{i}$ are uniformly bounded in $C^{m,\alpha}$ and 
$C^{m-2,\alpha}$ respectively, and invariant under rotations by an angle of 
order $r/R << 1$ on each circle of $T^{n-1}(r)$; here $r_{i}/R_{i} \rightarrow 0$ 
as $i \rightarrow \infty$. A function on a circle which is bounded in $C^{k}$ 
norm by 1, and which is periodic of period $\delta << 1$ is $\varepsilon$-close 
to its average value in $C^{k-1}$, where $\varepsilon$ depends linearly on 
$\delta$ for $\delta$ sufficiently small. 

 Moreover, in the region where $r(x) \sim r$, the black hole metric $g_{BH}$ 
differs from the cusp metric $g_{C}$ on the order of $O(r^{-(n-1)})$, cf. 
(3.8ff). It then follows that the equation (3.14), i.e. $(L(h))_{ab} = f_{ab}$, 
$a, b \geq 2$, may be written in the form
\begin{equation} \label{e3.31}
r^{2}(h_{ab}(r))'' + nr(h_{ab}(r))' = f_{ab}(r) + e_{ab}(r),
\end{equation}
where $e_{ab}(r) = O(r/R) + O(r^{-(n-1)})$, and as above the index $i$ 
has been supressed, (compare with (3.28)). 

  By (3.25), we already know that there exist $r_{i} \rightarrow \infty$ such 
that $r_{i}/R_{i} \rightarrow 0$ and $|h_{i}|(x) \rightarrow 0$ whenever 
$r(x) \geq r_{i}$. Hence view (3.31) for $r$ in the interval $[C_{0}, r_{i}]$, 
where $C_{0}$ is a fixed but arbitrarily large constant. The equation (3.31) 
may be integrated explicitly to give 
\begin{equation} \label{e3.32}
h_{ab}'(r) = \frac{1}{r^{n}}[\int_{C_{0}}^{r}r^{n-2}(f_{ab}(r) + e_{ab})dr + c_{1}],
\end{equation}
where we recall $h = h_{i}$, $f = f_{i}$. Let $\alpha_{i} = \sup f_{i}(r)\log R_{i}$ 
on the interval $[C_{0}, r_{i}]$, so that by (3.18), $\alpha_{i} \rightarrow 0$ as $i 
\rightarrow \infty$. Then (3.32) gives
$$|h_{ab}'(r)| \leq C[\frac{\alpha_{i}}{\log R_{i}}\frac{1}{r} + \frac{1}{R_{i}} + 
r^{-n}\log r],$$
on $[C_{0}, r_{i}]$. Integrating further from $r$ to $r_{i}$ then gives
\begin{equation} \label{e3.33}
|h_{ab}(r)| \leq C'[\alpha_{i}\frac{\log r_{i}}{\log R_{i}} + \frac{r_{i}}{R_{i}} + 
r_{i}^{-(n-1)}\log r_{i} + r^{-(n-1)}\log r] + |h_{ab}(r_{i})|
\end{equation}
$$\leq C'\delta_{i} + C'r^{-(n-1)}\log r + |h_{ab}(r_{i})|,$$
uniformly on $[C_{0}, r_{i}]$, where $\delta_{i} \rightarrow 0$ as $i \rightarrow 
\infty$. By (3.25), $|h_{i}(r_{i})| \rightarrow 0$ as $i \rightarrow \infty$. 

  This proves the claim (3.29), and as in Case (i), it follows that $h_{i}(x_{i}) 
\rightarrow 0$ as $i \rightarrow \infty$, for any $x_{i}$ satisfying (3.22). 

\medskip

  {\bf Case (iii).} 

   For $x_{i}$ satisfying (3.23), the metrics $(M_{i}, \widetilde g_{i}, x_{i})$ 
are also highly collapsed in regions of arbitrary but uniformly bounded diameter 
about $x_{i}$. However, just as above in Case (ii), the collapse may be unwrapped 
by passing to sufficiently large finite covering spaces and one may then pass to 
a limit. The limit is a complete black hole metric $g_{BH}$ on $D^{2}\times T^{n-2}$ 
as in (3.6). Similarly, as above, the forms $h_{i}$, (and $f_{i}$), lift to forms 
$\widetilde h_{i}$ on the covering spaces and converge, (in a subsequence), in the 
$\widetilde C^{m, \alpha'}$ topology, to a limit $T^{n-1}$-invariant form $\widetilde h$ 
satisfying the kernel equation (3.26) on $(D^{2}\times T^{n-2}, g_{BH})$. The 
assumption (3.18), together with the results above in Cases (i) and (ii) and 
the smooth convergence to the limit imply that one must have 
\begin{equation}\label{e3.34}
||\widetilde h||_{L^{\infty}} = 1 .
\end{equation}
In particular, $\widetilde h \neq 0$. Further, by (3.33) the limit form $\widetilde h$ 
satisfies
\begin{equation}\label{e3.35}
|\widetilde h| \leq C'r^{-(n-1)}\log r
\end{equation}
as $r \rightarrow \infty$ in $(D^{2}\times T^{n-2}, g_{BH})$. 

  The following Lemma now shows this situation is impossible.

\begin{lemma}\label{l3.4}
 Any bounded $T^{n-1}$-invariant Einstein deformation $h$ of a black 
hole metric $(D^{2}\times T^{n-2}, g_{BH})$ in (3.6) satisfies
\begin{equation}\label{e3.36}
|h|(y) \rightarrow  c_{0} \geq 0, \ \ {\rm as} \ \ y \rightarrow  
\infty, 
\end{equation}
for some constant $c_{0}$. Further, $c_{0} = 0$ if and only if $h 
\equiv 0$. In particular, the operator $L$ has trivial $L^{2}$ kernel, 
i.e. there are no non-trivial solutions $h$ of (3.26) with $h \in  
L^{2}$. 
\end{lemma}

{\bf Proof:} It is possible to prove Lemma 3.4 by a direct, although 
rather lengthy computation, by solving the system of ODE's for the coefficients 
of $h$ as in (3.28) above. Thus, the main point is to prove that $L$ has no 
$L^{2}$ kernel, i.e.~ the black hole metric is non-degenerate, cf. [26]. Since 
$g_{BH}$ has regions where the sectional curvature is positive when $n > 3$, 
this is not so easy to prove computationally. Thus, instead of going through the 
extensive computational details, we give a more conceptual proof at the non-linear 
level. 

  Thus, we first note that any complete Einstein metric (1.1) on 
$D^{2}\times T^{n-2}$ with an isometric $T^{n-1}$ action, with 
codimension 1 principal orbits, is a black hole metric $g_{BH}$ as in 
(3.6). This is proved in [5] when $n = 4$ and the same proof holds in 
all dimensions. A black hole metric is uniquely determined, up to isometry, 
by the flat structure induced on $T^{n-2}$, the mass parameter $m$, giving the 
length of the remaining $S^{1}$, (parametrized by $\theta$), and the homotopy 
class of $\sigma$. In particular, the only {\it small} deformations of $g_{BH}$ 
are those induced by variation of the flat structure on $T^{n-2}$ and variation 
of the mass $m$, cf. (2.6).

 Next we claim that the infinitesimal deformation $h$ is tangent to the 
moduli space of $C^{2}$ conformally compact (or asymptotically 
hyperbolic) Einstein metrics on the given manifold. To see this, since 
$h$ is invariant with respect to the standard $T^{n-1}$-action on $g_{BH}$, 
it may be written in the form (3.27), i.e. 
$$h = \sum h_{ab}(r)\theta^{a} \cdot \theta^{b},$$
where $\theta^{a}$ is the natural co-framing of $g_{BH}$, dual to $e_{a}$ 
as in (2.8). As noted in (2.9), the function $\rho  = r^{-1}$ is a smooth 
defining function, and gives a smooth compactification $\bar g_{BH} = 
\rho^{2}g_{BH}$ of $g_{BH}$. The associated compactification $\bar h = 
\rho^{2}h$ of $h$ satisfies $|\bar h|_{\bar g_{BH}} = |h|_{g_{BH}}$. 
Further, the equation (3.24) for an infinitesimal Einstein deformation 
may be reexpressed in terms of the compactified metric $\bar g_{BH}$ 
and $\bar h$, where it gives a system of ODE's for the functions $\bar 
h_{ab}(\rho)$. Since $g_{BH}$ is asymptotic to the hyperbolic cusp metric, 
it is easy to see that to leading order, the system (3.24) has the same form 
as that for the hyperbolic cusp metric, given in (A.8), (A.10) and (A.12). 
Hence a straightforward calculation for conformal changes of metric shows 
the coefficients $\bar h_{ab}(\rho)$ satisfy
$$\bar h_{ab}'' - \frac{n-2}{\rho}\bar h_{ab}' = o(1),$$
when $a, b \geq 2$. A similar expression holds for the coefficients $h_{1a}$. 
It follows by elementary integration that $\bar h$ extends $C^{2}$ up to the 
boundary at $\rho = 0$. This means that $h$ defines a tangent vector to the 
space of conformally compact Einstein metrics, as required. (A similar but 
much more elementary argument holds when $n = 3$, using the fact that 
infinitesimal Einstein deformations are infinitesimal hyperbolic deformations; 
we will not carry out the details). 

  Now the space of such $C^{2}$ conformally compact Einstein metrics is 
a smooth Banach manifold, and any tangent vector $h$ is tangent to a 
curve of conformally compact Einstein metrics, cf. [2], [4]. Since 
$h$ is $T^{n-1}$ invariant, it follows by the classification above that 
$h$ is tangent to the space of black hole metrics on $D^{2}\times 
T^{n-2}$. Thus, $h$ corresponds to an infinitesimal deformation of the flat 
structure on $T^{n-2}$ and the mass $m$.

  Because $h$ is $T^{n-1}$ invariant near infinity, it is now clear that 
$|h| \rightarrow  c_{0}$ at infinity, for some constant $c_{0}$. This gives 
(3.36). To prove the second statement, suppose $h$ is non-trivial, i.e. 
$h \neq 0$. If $h$ induces a non-trivial deformation of the $T^{n-2}$ 
factor, then it is clear from the form of $h$ above that $c_{0} \neq 0$. 
If instead the variation of the $T^{n-2}$ factor is trivial, consider 
the deformation of the mass $m$. This induces a variation of the length 
$\beta$ of the $S^{1}$ factor parametrized by $\theta$. Since $h \neq 
0$, the variation of $m$ is non-trivial. Now as noted following (2.9), 
$\beta$ is strictly monotone decreasing in $m$, and from (1.11), 
$\beta'(m) < 0$. Hence, the variation of the $S^{1}$ factor is 
non-trivial. This implies that $c_{0} \neq 0$, which completes the 
proof of Lemma 3.4. 
{\endproof}

  Combining the results obtained in Cases (i)-(iii) above, this now also 
completes the proof of (3.15). To prove the last statement in Proposition 
3.2,  (3.15) implies that ${\rm Ker} L = 0$ on ${\mathbb S}_{2}^{m,\alpha}$ 
in the $\widetilde C^{m,\alpha}$ norm. Since $L$ is essentially self-adjoint, 
and $M$ is assumed compact, standard Fredholm theory implies that 
that $L$ is surjective onto ${\mathbb S}_{2}^{m-2,\alpha}$ with the 
$\widetilde C^{m-2,\alpha}$ norm. Moreover, (3.15) then gives a bound 
$\Lambda$ on the norm of the inverse mapping $L^{-1}$ on these spaces. 
{\endproof}

\begin{remark}\label{r3.5}
{\rm With some further work, it should be possible to give a direct, effective 
proof of Proposition 3.2, avoiding the use of a contradiction. However, 
this requires understanding of the possible limit behaviors discussed 
above anyway, and carrying along effective estimates at each stage of 
the proof. We do not know of any proof that holds without addressing 
the structure of the possible limits. 

  A more explicit estimate of the constant $\Lambda$ would give more 
precise information on the set of Dehn fillings which carry Einstein metrics. }
\end{remark}

  Next, we observe that the proof of Proposition 3.2 also shows that 
the conclusion (3.15) holds for all smooth metrics sufficiently close 
to the approximate solution $\widetilde g$. More precisely, let 
$B_{\widetilde g}(\varepsilon)$ be the $\varepsilon$-ball about 
$\widetilde g$ in the $\widetilde C^{m, \alpha}$ topology on ${\mathbb M}$, 
cf. (2.22). 
\begin{corollary}\label{c3.6}
There exists $\varepsilon_{0} > 0$ such that (3.15) holds, for all metrics 
$g' \in B_{\widetilde g}(\varepsilon_{0})$, with again $\Lambda$ independent 
of $\bar \sigma$, (provided $\bar \sigma$ is sufficiently large). 
\end{corollary}

{\bf Proof:}
 The proof is exactly the same as that of Proposition 3.2. Briefly, if 
not, then there exists a sequence $(M_{i}, \widetilde g_{i})$, together 
with symmetric forms $h_{i}$ such that (3.18) holds, for some sequence 
of metrics $g_{i}' \in B_{\widetilde g_{i}}(\varepsilon_{i})$, with 
$\varepsilon_{i} \rightarrow 0$. However, the proof of Proposition 3.2 
applies just the same to this sequence, (as with the sequence 
$\widetilde g_{i}$ before), and gives the same contradiction.

{\endproof}

\medskip

  {\bf Step III. (Solution of the Nonlinear Problem).} 

  We are now in position to prove Theorem 1.1. This is done first in the 
case (3.13) where all the ends of $N$ are capped by Dehn filling, 

\medskip

{\bf Proof of Theorem 1.1. ($M_{\bar \sigma}$ compact)}

  Let $M = M_{\bar \sigma} = M(\sigma_{1}, ...,\sigma_{q})$ be obtained 
from $N$ by Dehn filling all the toral ends of $N$. Let $\widetilde g$ 
be the approximate Einstein metric on $M$ constructed in Step I. By (3.10)-(3.11), 
$\Phi(\widetilde g) = 0$ outside the glueing region $U = \cup U_{j}$. Write 
$M \setminus U = B \cup N_{T}$, where $B$ is the union of the black hole regions 
and $N_{T} \subset N$. 

  Let 
\begin{equation} \label{e3.37}
{\mathcal W} = \{f \in {\mathbb S}_{2}^{m-2,\alpha}: f(x) = 0, \forall x 
\ {\rm s.t.} \ dist_{\widetilde g}(x, B) \geq 2\}. 
\end{equation}
Note that ${\mathcal W}$ is closed in ${\mathbb S}_{2}^{m-2,\alpha}$ and so 
is a Banach subspace of ${\mathbb S}_{2}^{m-2,\alpha}$. Set $f_{0} = 
\Phi_{\widetilde g}(\widetilde g)$, and note that $f_{0} \in {\mathcal W}$. 
We let ${\mathcal W}_{\varepsilon} = {\mathcal W}\cap B_{f_{0}}(\varepsilon)$, 
where $B_{f_{0}}(\varepsilon)$ is the $\varepsilon$-ball about $f_{0}$ in 
${\mathbb S}_{2}^{m-2,\alpha}$, and set 
$${\mathcal U}_{\varepsilon} = \Phi^{-1}({\mathcal W}_{\varepsilon}),$$
so that 
\begin{equation} \label{e3.38}
\Phi_{0} = \Phi |_{{\mathcal U}_{{\varepsilon}}}: {\mathcal U}_{{\varepsilon}} 
\rightarrow {\mathcal W}_{{\varepsilon}}.
\end{equation}

  By Proposition 3.2 and Corollary 3.6, for $\varepsilon_{0}$ sufficiently small, 
every point in ${\mathcal W}_{{\varepsilon}_{0}}$ is a regular value of $\Phi$ and 
so of $\Phi_{0}$. Hence by the inverse function theorem, ${\mathcal U}_{{\varepsilon}_{0}}$ 
is a Banach submanifold of ${\mathbb M}^{m,\alpha}$, (of infinite codimension), and 
$\Phi_{0}$ is a local diffeomorphism onto ${\mathcal W}_{\varepsilon_{0}}$.
Of course the use here of Proposition 3.2 and Corollary 3.6 means that $\varepsilon_{0}$ 
might depend on $\bar \sigma$, via $R_{max}$. Further, by construction, 
$$\widetilde g \in {\mathcal U}_{\varepsilon_{0}}.$$ 

  We now consider the mapping $\Phi_{0}$ in place of $\Phi$. Being the restriction 
of a smooth map to a submanifold, $\Phi_{0}$ is of course still smooth. The 
linearization $L = D\Phi$ restricted to the tangent spaces 
$T_{g'}{\mathcal U}_{\varepsilon_{0}}$ of ${\mathcal U}_{{\varepsilon}_{0}}$, 
gives a linear mapping
\begin{equation} \label{e3.39}
L_{0}(h) = f ,
\end{equation}
from $h \in T_{g'}{\mathcal U}_{{\varepsilon}_{0}}$ to $f \in 
T_{\Phi_{0}(g')}{\mathcal W}_{{\varepsilon}_{0}}$. Observe that $f$ has restricted 
support on $M$; ${\rm supp}~f \subset {\rm supp}~\eta$. Of course in general $h$ does 
not have this form; one may well have ${\rm supp}~h = M$. 

  We now claim that Proposition 3.2, (and Corollary 3.6), can be improved when 
$\Phi$ is restricted to $\Phi_{0}$. 

\begin{proposition} \label{p3.7}
  Let $M = M_{\bar \sigma}$ be compact as in (3.13). Then there exist 
$\varepsilon_{0} > 0$ and $\Lambda < \infty$, independent of $\bar \sigma$, 
such that for any $g' \in {\mathcal U}_{\varepsilon_{0}}$ and $\bar \sigma$ 
sufficiently large, one has
\begin{equation} \label{e3.40}
||h||_{{\widetilde C}^{m,\alpha}} \leq \Lambda ||f||_{{\widetilde C}^{m-2,\alpha}} ,
\end{equation}
for $h$ and $f$ as in (3.39). Thus, $L_{0}$ is invertible on 
${\mathcal U}_{\varepsilon_{0}}$, and one has a uniform bound $\Lambda$ 
for the norm of $L_{0}^{-1}$, independent of $\bar \sigma$. 
\end{proposition}

{\bf Proof:}
Given the work above, this is now essentially an immediate consequence of the 
proof of Proposition 3.2. Thus, suppose first that $g' = \widetilde g$. The proof 
that (3.40) holds at $\widetilde g$ then follows exactly the proof of Proposition 3.2, 
with $f \in T_{f_{0}}{\mathcal W}$ and $h \in T_{\widetilde g}{\mathcal U}$ in place 
of the general $f$ and $h$ from before. The $\log R_{\max}$ term in the estimate (3.15) 
arises only because of the behavior in (3.32) in Case (ii). For $f \in T_{f_{0}}{\mathcal W}$, 
one has $f \equiv 0$ in this region and hence the same analysis following (3.32) shows 
that (3.30) holds. The proof of Cases (i) and (iii) holds without any changes. This 
proves that (3.40) holds at $\widetilde g$. The proof that it also holds for 
$g' \in {\mathcal U}_{\varepsilon_{0}}$, with $\varepsilon_{0}$ independent of 
$\bar \sigma$ for $\bar \sigma$ sufficiently large, is then exactly the same as 
Corollary 3.6, with $\Phi_{0}$ in place of $\Phi$. The last statement also follows 
as before, since $L_{0}$ is still essentially self-adjoint as a mapping 
$T({\mathcal U}_{\varepsilon_{0}}) \rightarrow T({\mathcal W}_{\varepsilon_{0}})$. 
{\endproof}

  The proof of Theorem 1.1 when $M$ is compact is now quite straightforward. 
First, the estimate (3.11) implies that
\begin{equation} \label{e3.41}
||\Phi_{0}(\widetilde g)||_{\widetilde C^{m-2, \alpha}} \leq (R_{min})^{-(n-1)},
\end{equation}
where via (3.4), $R_{min}$ is the shortest length of the collection of 
geodesics $\{\sigma_{j}\}$ in $\bar \sigma$, up to a fixed constant. 

  Next, let $h$ be any symmetric bilinear form in $T_{\widetilde g}{\mathcal U}$ 
satisfying $||h||_{\widetilde C^{m, \alpha}} \leq 1$. Then (3.11) and (3.12) show 
that
$$||D_{\widetilde g}\Phi_{0}(h)||_{\widetilde C^{m-2, \alpha}} \leq K.$$
The constant $K$ depends only on the local geometry of $\widetilde g$, 
(in covering spaces in sufficiently collapsed regions), and hence is 
independent of $\bar \sigma$. For the same reasons, choosing $\varepsilon_{0} > 0$ 
smaller if necessary, one has
\begin{equation}\label{e3.42}
||D_{g'}\Phi_{0}(h)||_{\widetilde C^{m-2, \alpha}} \leq 2K,
\end{equation}
for all $g' \in {\mathcal U}_{\varepsilon_{0}}$, and $h$ as above, where 
$K$ is independent of $\bar \sigma$. 
Next, Proposition 3.7 shows that 
\begin{equation}\label{e3.43}
||(D_{g'}\Phi_{0})^{-1}(f)||_{\widetilde C^{m, \alpha}} \leq \Lambda ,
\end{equation}
for all $g' \in {\mathcal U}_{{\varepsilon}_{0}}$ and $f \in 
T{\mathcal W}_{\varepsilon_{0}}$ with $||f||_{\widetilde C^{m-2, \alpha}} \leq 1$. 
The bounds (3.42)-(3.43) prove that $\Phi_{0}$ is a bi-Lipschitz map, with Lipschitz 
constant $2K$ for $\Phi_{0}$ and $\Lambda$ for $\Phi_{0}^{-1}$. 

  The inverse function theorem applied to the mapping $\Phi_{0}$ between 
the Banach manifolds ${\mathcal U}_{{\varepsilon}_{0}}$ and 
${\mathcal W}_{{\varepsilon}_{0}}$ then implies that there is a domain $\Omega 
\subset U_{\varepsilon_{0}}$ and $\varepsilon_{1} > 0$ such that 
\begin{equation}\label{e3.44}
\Phi_{0}: \Omega \rightarrow {\mathcal W}_{\varepsilon_{1}} ,
\end{equation}
is a diffeomorphism onto ${\mathcal W}_{\varepsilon_{1}}$. The constant 
$\varepsilon_{1}$ is of the form $\varepsilon_{1} = (4K/\Lambda)\varepsilon_{0}$. 
By (3.41), one may now choose $R_{min}$ sufficiently large, i.e. $\bar \sigma$ 
sufficiently large, so that $0 \in {\mathcal W}_{\varepsilon_{1}}$. 
Via (3.44), this implies that there exists a metric 
$g \in {\mathcal U}_{\varepsilon_{0}}$, such that
$$\Phi(g) = \Phi_{0}(g) = 0.$$

  By Lemma 2.1, $g$ is then an Einstein metric on $M$, smoothly close to 
$\widetilde g$. 
{\endproof}

\begin{remark}\label{r3.8}
  {\rm Since $\Phi_{0}$ in (3.44) is a diffeomorphism on 
$\Omega$, the metric $g$ is the unique Einstein metric, (up to isometry), 
with the normalization (1.1) in $\Omega$. Moreover, since $\Phi$ is a local 
diffeomorphism near $g$, it follows that the metrics $g$ constructed above 
are isolated points in the moduli space of Einstein metrics on $M$. }
\end{remark}

\medskip

  Next, we complete the proof of Theorem 1.1 by discussing the case where 
not all the cusps of $M$ are capped by Dehn filling. 

\medskip

{\bf Proof of Theorem 1.1. ($M_{\bar \sigma}$ non compact)}

  Let ${\mathcal E}_{c} ={\mathcal E}_{c}(N)$ be the collection of Einstein metrics 
constructed on the compact manifolds $M_{\bar \sigma}$ above associated to a given 
$N$. This is an infinite collection of metrics, parametrized by $\bar \sigma$. 
Now let $M = M_{\bar \sigma} = M(\sigma_{j_{1}}, ..., \sigma_{j_{p}})$ be any manifold 
obtained by Dehn filling a collection of $p$ toral ends $E_{j_{k}}$ of $N$, with each 
$\sigma_{j_{k}}$ sufficiently large. By relabeling, assume $1 \leq j_{k} \leq p$, 
so that the ends $E_{j}$, $p+1 \leq j \leq q$ are cusp ends of $M$. Further, we 
assume $p < q$, so that $M$ is non-compact. The manifold $M$ may be written in the 
form $M = M(\sigma_{1}, ..., \sigma_{p}, \infty, ..., \infty)$. 

  Let $M_{i} = M(\sigma_{1}, ..., \sigma_{p}, \sigma_{p+1}^{i}, \cdots, \sigma_{q}^{i})$, 
where $\sigma_{k}^{i}$, $p+1 \leq k \leq q$, is any sequence such that 
$\sigma_{k}^{i} \rightarrow \infty$ as $i \rightarrow \infty$, for each fixed $k$. 
Let $\widetilde g_{i}$ be the approximate Einstein metrics constructed on $M_{i}$ and 
let $g_{i} \in {\mathcal E}_{c}$ be the associated Einstein metrics on $M_{i}$ given 
by Theorem 1.1, (in the compact case). If $y_{0}$ is any fixed point in $N$, it is 
clear that the pointed sequence $(M_{i}, \widetilde g_{i}, y_{0})$ has a subsequence 
converging smoothly and uniformly on compact sets to the limit manifold 
$(M, \widetilde g, y_{0})$, where $\widetilde g$ is the approximate Einstein metric 
constructed on $M$ in Step I. Since the Einstein metrics $g_{i}$ are smoothly close 
to the approximate metrics $\widetilde g_{i}$, $\{g_{i}\}$ also converges, again smoothly 
and uniformly on compact sets, to a limit Einstein metric $g$ on $M$. The limit 
$g$ is complete, and of uniformly bounded curvature. This completes the proof of 
Theorem 1.1. 
{\endproof}

  Having completed the proof of Theorem 1.1, we next show that the 
homeomorphism type of the Dehn-filled manifolds $M_{\bar \sigma}$ is 
determined up to finite ambiguity by the data $\bar \sigma = 
(\sigma_{1}, ..., \sigma_{p})$. Let Out$(\pi_{1}(N))$ be the group of 
outer automorphisms of $\pi_{1}(N)$. By Mostow-Prasad rigidity, this is 
a finite group, isomorphic to the isometry group Isom$(N, g_{-1})$ of 
$N$.

\begin{proposition}\label{p3.9}
Let $n \geq 4$. The number of manifolds $M_{\bar \sigma}$ homeomorphic 
to a given manifold $M_{\bar \sigma_{0}}$ is finite, and 
bounded by the cardinality of Out($\pi_{1}(N))$.
\end{proposition}

{\bf Proof:} If $M = M_{\bar \sigma}$ is obtained from $N$ by Dehn 
filling a collection of cusp ends $\{E_{j}\}$ of $N$, then by the 
Seifert-Van Kampen theorem, the fundamental group $\pi_{1}(M)$ is given 
by 
\begin{equation}\label{e3.45}
\pi_{1}(M) = \pi_{1}(N)/\langle \cup R_{j} \rangle,
\end{equation}
where $R_{j} \simeq {\mathbb Z}$ is the subgroup generated by the 
closed geodesic $\sigma_{j} \in \bar \sigma$, (i.e. the meridian circle 
is annihilated). As noted in \S 2.1, if the Dehn filling is 
sufficiently large, then $M$ is a $K(\pi, 1)$ and each core torus 
injects in $\pi_{1}$:
$$\pi_{1}(T_{j}^{n-2}) \hookrightarrow \pi_{1}(M).$$

  Thus to each peripheral subgroup ${\mathbb Z}^{n-1} \simeq 
\pi_{1}(E_{j}) \subset \pi_{1}(N)$ is associated a subgroup ${\mathbb 
Z}^{n-2} \subset \pi_{1}(M)$, obtained by dividing ${\mathbb Z}^{n-1}$ 
by ${\mathbb Z}$. This gives a distinguished collection of (conjugacy 
classes of) subgroups isomorphic to ${\mathbb Z}^{n-2}$ and 
${\mathbb Z}^{n-1}$, corresponding to the filled and unfilled ends 
of $N$; call these the peripheral subgroups of $\pi_{1}(M)$. As before 
with $N$, any non-cyclic abelian subgroup of $\pi_{1}(M)$ is conjugate to 
a subgroup of a peripheral subgroup. This is because $M$ admits a complete 
metric of non-positive sectional curvature naturally associated to the Dehn 
filling, cf. \S 2.1. With respect to such a metric, any non-cyclic 
abelian subgroup is carried by an essential torus embedded in $M$. 
However, up to isotopy, all such tori are contained in the core tori 
$T^{n-2}$ of $M$ or  the end tori $T^{n-1}$ of $M$.

  Now suppose $M_{i}$, $i = 1,2$, are two $n$-manifolds obtained by 
Dehn fillings of a given hyperbolic $N$. If $M_{1}$ is homeomorphic 
to $M_{2}$, then $\pi_{1}(M_{1}) \simeq \pi_{1}(M_{2})$, and 
we may choose a fixed isomorphism identifying both with the (abstract) 
group $\pi_{1}(M)$. A homeomorphism $F: M_{1} \rightarrow M_{2}$ 
then defines an automorphism 
\begin{equation}\label{e3.46}
F_{*}: \pi_{1}(M) \rightarrow \pi_{1}(M). 
\end{equation}

  By the uniqueness mentioned above, it follows that $F_{*}$ permutes 
the collection of peripheral subgroups onto themselves, inducing an 
isomorphism of each ${\mathbb Z}_{i}^{n-2}$ to some ${\mathbb 
Z}_{j}^{n-2}$, and ${\mathbb Z}_{k}^{n-1}$ to some ${\mathbb Z}_{l}^{n-1}$ 
up to conjugacy; of course one may have $i = j$ or $k = l$. Each such 
subgroup is carried by an embedded, essential torus $T^{n-2}$ or 
$T^{n-1}$ in $M$. Let $\hat T_{i}^{n-2} = F(T_{i}^{n-2})$ and set 
$T = \cup T_{i}^{n-2}$, $\hat T = \cup \hat T_{i}^{n-2}$. Then $F$ gives 
a homeomorphism of $N = M \setminus T$ onto $\hat N = M \setminus \hat T$. 
Equivalently, $F$ induces a homeomorphism of the original hyperbolic 
manifold $N$, 
\begin{equation}\label{e3.47}
F: N \rightarrow N,
\end{equation}
permuting the cusp ends of $N$. Further, if $F$ maps the end $E_{i}$ to 
$E_{j}$ then by (3.46), $F_{*}\langle \sigma_{i} \rangle = \langle 
\sigma_{j} \rangle$, up to conjugacy, in $\pi_{1}(N)$; here $\langle 
\sigma \rangle$ is the subgroup generated by $[\sigma]$. 

  If $F$ is homotopic to the identity on $N$, then the filling data of 
$M_{1}$ and $M_{2}$ are the same, up to sign, and so $M_{1}$ and 
$M_{2}$ are diffeomorphic, cf. \S 2.1. If not, then $F$ induces a 
non-trivial automorphism $F_{*}$ of $\pi_{1}(N)$, so that $F_{*}$ is an 
element of the outer automorphism group Out$(\pi_{1}(N))$. Since this 
group is finite, it follows that only a finite number of filling data 
can give rise to homeomorphic manifolds $M_{\bar \sigma}$. One obtains a 
bound on this number by a bound on the order of Isom$(N)$, or more 
precisely a bound on the order of the corresponding effective group 
acting on the corresponding Dehn filling spaces ${\mathbb Z}^{n-1}$. 
{\endproof}

\bigskip

  We complete this section with a discussion of Dehn filling on 
non-toral ends. Thus, let $(N, g_{-1})$ be a complete hyperbolic 
$n$-manifold of finite volume, with an end $E$ of the form $F\times 
{\mathbb R}^{+}$, where $F$ is a flat manifold with induced metric 
$g_{0}$. By the Bieberbach theorem, cf. [34],
\begin{equation} \label{e3.48}
F = T^{n-1} / \Gamma,
\end{equation}
where $\Gamma$ is a finite group of Euclidean isometries acting freely 
on $T^{n-1}$. Let $\bar E$ be the covering space of $E$ with covering 
group $\Gamma$, so that $\bar E$ is of the form $T^{n-1}\times {\mathbb 
R}^{+}$, with hyperbolic metric $g_{-1}$. For $\sigma$ a simple closed 
geodesic in $(T^{n-1}, g_{0})$, let $\phi_{\sigma}$ be a diffeomorphism 
of $\partial(D^{2}\times T^{n-2})$ to $T^{n-1}$ sending $S^{1} = 
\partial D^{2}$ to $\sigma \subset T^{n-1}$, so that $\phi_{\sigma}$ 
attaches a solid torus to $T^{n-1}$ along $\sigma$. Now suppose that 
the action of $\Gamma$ on $T^{n-1}$ extends to a free action of 
$\Gamma$ on $D^{2}\times T^{n-2}$ and that $\Gamma$ commutes with the 
diffeomorphism $\phi_{\sigma}$ on the boundary $T^{n-1}$. Then the 
quotient manifold
$$M_{\sigma} = (D^{2}\times T^{n-2})/\Gamma \cup_{\phi_{\sigma}} N$$
is well-defined, and is the manifold obtained by performing Dehn 
filling the end $E$ along the geodesic $\pi(\sigma) \subset F$, where 
$\pi: T^{n-1} \rightarrow F$ is the covering projection. 

  The following result gives a necessary and sufficient condition for 
the existence of such Dehn fillings of an end $E$. 

\begin{lemma}\label{l3.10}
For $F$ and $\sigma$ as above, the quotient $M_{\sigma}$ is 
well-defined, and carries a corresponding quotient of the AdS black 
hole metric $g_{BH}$ in (3.6) if and only if, for any $\gamma \in 
\Gamma$ acting on the universal cover ${\mathbb R}^{n-1}$, one has
\begin{equation}\label{e3.49}
\langle \gamma(\sigma) \rangle \parallel \langle \sigma \rangle,
\end{equation}
where $\langle \tau \rangle$ is the line through $\tau$.
\end{lemma}

{\bf Proof:} In the process of Dehn filling a toral end, the initial 
flat structure on $T^{n-1}$ is deformed along a curve of flat 
structures, by smoothly changing the length of the meridian curve 
$\sigma$ from its initial length to length 0. This is described 
explicitly in (3.7). Thus, one has to check if the deformation (3.7) is 
invariant under a corresponding deformation of the action of $\Gamma$.

  As discussed following (3.5), let $(\sigma, b_{2}, ..., b_{n-1})$ be 
a basis for the lattice giving $T^{n-1}$, and set $\sigma = b_{1}$. Let 
$b_{i}^{r} = b_{i}+(\lambda(r)-1)(\langle b_{i},\sigma \rangle / 
|\sigma|^{2})\sigma$ be as in (3.7), and let $t_{b_{i}}^{r}$ denote the 
generators for the lattice $({\mathbb Z}^{n-1})(r)$ defining 
$T^{n-1}(r)$; thus $t_{b_{i}}^{r}$ is translation by the vector 
$b_{i}(r)$ on ${\mathbb R}^{n-1}$.  

  By the Bieberbach theorem (3.48), the group $\pi_{1}(F)$ is a 
semi-direct product of ${\mathbb Z}^{n-1}$ with $\Gamma$. The group 
$\Gamma$ acts by affine transformations on ${\mathbb R}^{n-1}$; each 
$\gamma \in \Gamma$ acts by $(A_{\gamma}, t_{\gamma})$, where 
$A_{\gamma} \in O(n-1)$ and $t_{\gamma}$ is a translation on ${\mathbb 
R}^{n-1}$ by the vector $t_{\gamma}$. Thus $\gamma(v) = A_{\gamma}(v) + 
t_{\gamma}$ and
\begin{equation} \label{e3.50}
(A_{\gamma_{1}}, t_{\gamma_{1}})(A_{\gamma_{2}}, t_{\gamma_{2}})(v) = 
A_{\gamma_{1}}A_{\gamma_{2}}(v) + A_{\gamma_{1}}(t_{\gamma_{2}}) + 
t_{\gamma_{1}} = (A_{\gamma_{1}}A_{\gamma_{2}}, 
t_{A_{\gamma_{1}}(t_{\gamma_{2}})+ t_{\gamma_{1}}})(v).
\end{equation}
Define then a deformation of the action of $\Gamma$ by setting
\begin{equation} \label{e3.51}
A_{\gamma}^{r} = A_{\gamma}, \ {\rm and} \ t_{\gamma}^{r} = t_{\gamma} 
+ (\lambda(r)-1)\frac{\langle t_{\gamma}, \sigma 
\rangle}{|\sigma|^{2}}\sigma = t_{\gamma}^{\perp} + 
\lambda(r)\frac{\langle t_{\gamma}, \sigma \rangle}{|\sigma|^{2}}\sigma,
\end{equation}
where $t_{\gamma}^{\perp}$ is the component of $t_{\gamma}$ orthogonal 
to $\langle \sigma \rangle$. Thus, the orthogonal part $A_{\gamma}$ of 
$\gamma$ remains unchanged, while the translation part $t_{\gamma}^{r}$ 
varies along $\sigma$, and is orthogonal to $\sigma$ at $r = r_{+}$, 
where $\lambda(r_{+}) = 0$. Observe that the deformation 
$t_{\gamma}^{r}$ has exactly the same form as $t_{b_{i}}^{r}$.

  To verify that this gives a well-defined action of $\pi_{1}(F)$ on 
${\mathbb R}^{n-1}$ one needs to check that the relations of 
$\pi_{1}(F)$ are preserved. This is clear for the orthogonal (or $A$) 
part of the action by (3.50)-(3.51), and so one only needs to consider 
the translation or vector part of the action.

  Each relation $R$ is a word in some generators $A_{\gamma}$, 
$t_{\gamma}$, $t_{b_{i}}$. Thus, as a vector, $R(A_{\gamma}, 
t_{\gamma}, t_{b_{i}}) = 0$, where each $t$ acts by translation, (i.e. 
addition), and each $A_{\gamma}$ acts by an orthogonal matrix on some 
$t$ vector. To verify that $R^{r} = R(A_{\gamma}, t_{\gamma}^{r}, 
t_{b_{i}}^{r}) = 0$, suppose first that $R$ involves no rotational 
part, i.e. $R = R(t_{\gamma}, t_{b_{i}}) = 0$. The components of $R$ 
parallel and orthogonal to $\sigma$ then also both vanish. Since the 
deformations $t_{\gamma}^{r}$ and $t_{b_{i}}^{r}$ have exactly the same 
form along these components, and orthogonal projection commutes with 
translation, it follows that $R^{r} = R(t_{\gamma}^{r}, t_{b_{i}}^{r}) 
= 0$. 

  Next, consider the action of any $A = A_{\gamma}$ on some translation 
$t = t_{\gamma}$ or $t_{b_{i}}$. The condition (3.49) implies that $A$ 
leaves the subspaces $\langle \sigma \rangle$ and $\langle \sigma 
\rangle^{\perp}$ invariant, i.e. $A(\sigma) = \pm \sigma$. As above, 
the components of the vector $R = R(A_{\gamma}, t_{\gamma}, t_{b_{i}})$ 
along $\langle \sigma \rangle$ and $\langle \sigma \rangle^{\perp}$ 
vanish. Since any $A$ commutes with translation by $\sigma$, it follows 
that $R_{\sigma} = R((A_{\gamma})_{\sigma}, (t_{\gamma})_{\sigma}, 
(t_{b_{i}})_{\sigma}) = 0$, where $t_{\sigma}$ is the $\sigma$ 
component of $t$ and $A_{\sigma} = A|_{\langle \sigma \rangle}$. The 
same statement holds with respect to $\langle \sigma \rangle^{\perp}$. 
Since, as above, the vectors $t_{\gamma}^{r}$ and $t_{b_{i}}^{r}$ have the 
same form, it follows that the $\sigma$ and $\sigma^{\perp}$ components of 
$R^{r}$ also vanish, as required. This shows that the condition (3.49) 
is a sufficient condition that $M_{\sigma}$ is well-defined.

  Observe that the action of $\Gamma$ is well-defined at the core 
$(n-2)$-torus $T^{n-2} = \{r = r_{+}\}$ where $\lambda(r_{+}) = 0$, and 
so
\begin{equation}\label{e3.52}
\langle \gamma(b_{i}(r_{+})), \sigma \rangle = 0.
\end{equation}
Conversely, the condition (3.52) is necessary for the Dehn filling 
$M_{\sigma}$ to be well-defined. Since $\Gamma$ acts by isometries, 
$\langle \gamma(b_{i}(r_{+})), \sigma \rangle = \langle b_{i}(r_{+}), 
\gamma^{-1}\sigma \rangle$. However, by construction, i.e. (3.7), we 
know that $\langle b_{i}(r_{+}), \sigma \rangle = 0$, $\forall i > 1$. 
Hence, (3.52) requires the condition (3.49), so that (3.49) is also 
necessary. 
{\endproof}

 Define the Dehn filling along $\sigma$ to be {\it admissible} if 
$\Gamma$ and $\sigma$ satisfy the condition (3.49). This leads to the 
following extension of Theorem 1.1. 

\begin{corollary}\label{c3.11}
Let $M_{\bar \sigma}$ be any manifold obtained by performing a 
sufficiently large, admissible Dehn filling of the ends $E_{j}$, 
$1 \leq j \leq q$, of a complete hyperbolic $(N, g_{-1})$. 
Then $M_{\bar \sigma}$ admits an Einstein metric $g$ satisfying (1.1). 
\end{corollary}

{\bf Proof:} Using Lemma 3.10, one constructs the approximate Einstein 
metric $\widetilde g$ exactly as in Proposition 3.1. The rest of the 
proof proceeds exactly as in the proof of Theorem 1.1.
{\endproof}

  For a given end $E = F \times {\mathbb R}^{+}$ with $F = 
T^{n-1}/\Gamma$, not all Dehn fillings will be admissible, unless $E$ 
is toral. Nevertheless, for many such $F$, there will be an infinite 
number of admissible fillings; this can be checked by inspection.

\section{Further Results and Remarks}
\setcounter{equation}{0}

  In this section, we collect a number of remarks on the geometry and 
topology of the Einstein metrics $(M_{\bar \sigma}, g)$ constructed in 
Theorem 1.1 or Corollary 3.11, and prove the remaining results stated 
in the Introduction; Theorem 1.2 is proved in \S 4.1.

\medskip

{\bf \S 4.1.} 
By the Chern-Gauss-Bonnet theorem [15], if $N$ is a complete hyperbolic 
$n$-manifold of finite volume, then
\begin{equation}\label{e4.1}
vol N = (-4\pi)^{m}\frac{m!}{(2m)!}\chi(N),
\end{equation}
where $n = 2m$ and $\chi(N)$ is the Euler characteristic of $N$. In 
particular, the sign of the Euler characteristic is $(-1)^{m}$. Since 
the Dehn-filled manifold $M = M_{\bar \sigma}$ decomposes as a union of 
$N$ and a collection of solid tori $D^{2}\times T^{n-1}$, an elementary 
Mayer-Vietoris argument shows that 
$$\chi(N) = \chi(M).$$ 
Since $\chi(N)$ can be arbitrarily large for hyperbolic manifolds, (by 
passing to covering spaces), $\chi(M)$ can thus be made arbitrarily 
large when $n$ is even. 

\medskip

  Next we verify the claims (1.4) and (1.5). Regarding (1.4), the 
curvature of the black hole metric is given by (2.8), while that of the 
approximate Einstein metric $\widetilde g$ is as stated in Proposition 
3.1. The Einstein metric $g$ on $M$ is close to $\widetilde g$ in the 
$\widetilde C^{m,\alpha}$ topology, for any $m$. Hence, the curvature 
of $g$ is uniformly close to that of $\widetilde g$. This gives the 
estimate (1.4). 

   Regarding the Weyl curvature estimate (1.5), $W = 0$ on any 
hyperbolic manifold. For the black hole metric, as noted following 
(2.8), $W$ decays as $|W| = O(e^{-(n-1)s})$, where $s$ is the distance 
to the core $T^{n-2}$. On the other hand, the volume of the region 
$D(s)$ with respect to the approximate solution $\widetilde g$ is on the 
order of $O(e^{(n-1)(s-\ln R)})$, where $R$ is given by (3.4). It follows 
that the volume of the region where $|W| \geq \delta$ is on the order 
of $R^{-(n-1)}\delta^{-1}$. This verifies (1.5) for the approximate 
solution $\widetilde g$. Again, since the Einstein metric $g$ is 
uniformly close to $\widetilde g$, (1.5) follows for $g$. On the other 
hand, there is a fixed constant $c_{0} > 0$, depending only on 
dimension, such that
\begin{equation}\label{e4.2}
|W|_{L^{\infty}} \geq c_{0},
\end{equation}
since this is the case for the black hole metric $g_{BH}$ near the core 
torus $T^{n-2}$. Of course (4.2) assumes $n \geq 4$. 

  An immediate consequence of (1.5) and the Chern-Weil theory is that 
all Pontryagin numbers of $M$ vanish when $M$ is compact. In 
particular, by the Hirzebruch signature theorem, the signature $\tau(M) 
= 0$.

\begin{remark}\label{r4.1}
{\rm In a natural sense, most of the Einstein manifolds constructed are 
not locally isometric. (All hyperbolic manifolds are of course locally 
isometric). Let $N$ be a complete, noncompact hyperbolic manifold of 
finite volume, and let $\bar N$ be a covering of $N$ of degree $k$. If 
$M_{\bar \sigma}$ is obtained from $N$ by Dehn filling, then $M_{\bar 
\sigma}$ admits a degree $k$ covering $\bar M_{\bar \sigma}$, such that 
$\bar M_{\bar \sigma}$ is obtained from $\bar N$ by Dehn filling on 
cusps of $\bar N$; these Dehn fillings are lifts of the Dehn fillings 
on $M_{\bar \sigma}$. However, $\bar N$ admits many new Dehn fillings 
which are not lifts of Dehn fillings on $N$. Hence, ``most all'' of the 
Einstein metrics associated with $\bar N$ are not lifts of Einstein 
metrics associated to $N$. }
\end{remark}

\begin{remark}\label{r4.2}
{\rm Let $N$ be as above, and suppose $\pi_{1}(N)$ admits a 
homomorphism onto a free group $F_{2}$ with two generators. The lower 
bound in (1.3) is achieved by taking coverings of hyperbolic manifolds 
which admit such a surjection onto $F_{2}$, cf. [13], [24]. Let $C(N)$ 
denote the number of cusps of $N$. We claim that many coverings $\bar 
N^{k}$ of $N$ of degree $k$ have $C(\bar N^{k})$ growing linearly with 
$k$, i.e. linearly in the volume. More precisely, there exist 
constants, $c, d > 0$, depending only on dimension $n$, such that 
\begin{equation}\label{e4.3}
C(\bar N^{k}) \geq d\cdot k,
\end{equation}
for a collection of isometrically distinct coverings $\bar N^{k}$ of 
cardinality at least $e^{ck\ln k}$. 
  To see this, let $\phi: \pi_{1}(N) \rightarrow F_{2}$ be the 
surjective homomorphism onto $F_{2}$. Any subgroup $H$ of index $k$ in 
$F_{2}$ determines a covering space $\bar N^{k}$, with $\pi_{1}(\bar 
N^{k}) = (\phi)^{-1}(H)$. Since $F_{2}$ is free, $\phi$ sends any 
$\pi_{1}(T_{j}^{n-1}) \simeq {\mathbb Z}^{n-1}$ to $\langle a_{j} 
\rangle$, for some fixed $a_{j} \in F_{2}$. If $a_{j} \in H$, then the 
covering $\bar N^{k}$ unwraps $T_{j}^{n-1}$ into $k$ disjoint copies of 
$T^{n-1}$, giving rise to $k$ cusp ends, and thus giving (4.3). Hence, 
one needs to count the number of distinct index $k$ subgroups of 
$F_{2}$ containing a given element $a$. Following [20], there are at 
least $k\cdot k!$ subgroups of $F_{2}$ of index $k$, and at least $k!$ 
of these contain a given element $a \in F_{2}$. Following [13], this 
gives the lower bound on $c$ above for the number of non-isometric 
coverings.

   The opposite bound to (4.3),
$$C(\bar N^{k}) \leq D\cdot k,$$
for some fixed constant $D= D(n)$, is an immediate consequence of the 
Margulis Lemma. }
\end{remark}

  Next we prove the following expanded version of Theorem 1.2. 
Let ${\mathcal E}$ denote the class of Einstein metrics constructed via 
Theorem 1.1 or Corollary 3.11, together with the class of complete, 
non-compact hyperbolic $n$-manifolds $(N, g_{-1})$ of finite volume. 
\begin{theorem}\label{t4.3}
The space ${\mathcal E}$ is closed with respect to the pointed Gromov-Hausdorff 
topology or the $C^{\infty}$ topology. Further the volume functional 
\begin{equation}\label{e4.4}
vol: {\mathcal E} \rightarrow {\mathbb R}^{+}
\end{equation}
is continuous and proper with respect to these topologies. Any limit 
point $(M, g)$ of a sequence $(M^{i}, g^{i}) \in {\mathcal E}$ is 
obtained by opening a finite number of cusps of $M_{i}$. 
\end{theorem}

{\bf Proof:}  Let $(M^{i}, g^{i})$ be any sequence in ${\mathcal E}$ of 
bounded volume. By passing to a subsequence, we may assume that $M^{i} 
= M_{\bar \sigma^{i}}$, where $M_{\bar \sigma^{i}}$ is obtained from a 
fixed complete hyperbolic manifold $N$ by Dehn filling a collection of 
cusp ends. The sequence $\bar \sigma_{j}$ (partially) diverges to infinity 
in the Dehn filling space; thus for one and possibly several fixed $j$, 
$L(\sigma_{j}^{i}) \rightarrow \infty$ as $i \rightarrow \infty$, where 
$\sigma_{j}^{i}$ is a sequence of simple closed geodesics in tori 
$T_{j}^{n-1}$ in the $j^{\rm th}$ cusp end of $N$. By passing to a further 
subsequence, we may then assume that $M^{i}$ is obtained by Dehn filling 
of $a+b$ fixed cusps of $N$, and that $L(\sigma_{j}^{i}) \rightarrow \infty$ 
for $1 \leq j \leq a$, while $L(\sigma_{j}^{i})$ remains bounded, for 
$a+1 \leq j \leq a+b$. Here $a+b \leq q$, where $q$ is the number of cusps 
of $N$. 

  By construction, each Einstein metric $g^{i} \in B_{\widetilde g_{i}}(\varepsilon_{0})$, 
where $B(\varepsilon_{0})$ is the $\varepsilon_{0}$-ball in the $\widetilde C^{m,\alpha}$ 
topology and $\widetilde g_{i}$ is the approximate metric constructed on $M^{i}$; see 
the proof of Theorem 1.1. Further, as in the proof of Theorem 1.1 in the 
non-compact case, the sequence of metrics $\widetilde g_{i}$ 
converges, in a subsequence, to a limit metric $\widetilde g_{\infty}$ on a 
manifold $M_{\infty} = M_{\bar \sigma_{\infty}}$, where $\bar \sigma_{\infty} = 
(\infty, ..., \infty, \sigma_{a+1}, ..., \sigma_{a+b})$. Thus, $M_{\infty}$ 
is obtained from $M^{i}$ by opening $a$ cusps. The metric $\widetilde g_{\infty}$ 
is thus obtained from $N$ by Dehn filling $b$ cusps of $N$, along the curves 
$\sigma_{a+1}, ..., \sigma_{a+b}$. 

  Theorem 1.1, (or Corollary 3.11), thus gives the existence of an Einstein metric 
$g_{\infty}$ on $M_{\infty}$, in $B_{\widetilde g_{\infty}}(\varepsilon_{0})$. 
This proves that $\mathcal E$ is closed in the pointed $\widetilde C^{m,\alpha}$ 
topology and in fact ${\mathcal E}_{V} = \{g \in {\mathcal E}: vol_{g}M \leq V\}$ 
is compact. The convergence in the $C^{\infty}$ topology then follows from 
well-known elliptic regularity associated to the Einstein equation. The 
$C^{\infty}$ topology is much stronger than the Gromov-Hausdorff topology, 
hence ${\mathcal E}$ is also closed in the Gromov-Hausdorff topology. 

  To see that the volume functional (4.4) is continuous, the sequence 
$(M^{i}, g^{i})$ or $(M^{i}, \widetilde g^{i})$ converges smoothly to 
its limit, uniformly on compact subsets. Hence, for any compact domain 
$D \subset M_{\infty}$, $vol_{g^{i}}D \rightarrow vol_{g_{\infty}}D$. 
Further, if $D$ contains a sufficiently large region of $N$, the volume 
of the complement is uniformly small, for all $i$; this follows since 
the volume of the approximate metrics $\widetilde g$ at geodesic 
distance $t$ from the glueing tori is uniformly exponentially small. 
This proves the continuity of $vol$ on ${\mathcal E}$. The 
properness of $vol$ follows from the argument above: any sequence in 
${\mathcal E}$ of bounded volume has a convergent subsequence in 
${\mathcal E}$. Similarly, the fact that limits are obtained by opening 
cusps has already been proved above. 
{\endproof}

\bigskip

{\bf \S 4.2.}  In this section, we discuss further aspects of the volume and 
convergence behavior of the Einstein metrics constructed above in dimension 4. 
To begin, the Chern-Gauss-Bonnet theorem in dimension 4 states
\begin{equation}\label{e4.5}
\chi(M) = \frac{1}{8\pi^{2}}\int_{M}(|W|^{2} - \tfrac{1}{2}|z|^{2} + 
\tfrac{1}{24}s^{2})dV,
\end{equation}
where $z = Ric - \frac{s}{4}g$ is the trace-free Ricci curvature. The 
formula (4.5)  holds for all compact manifolds $M$. It also holds for 
complete non-compact hyperbolic manifolds of finite volume. This follows 
by using the Chern-Gauss-Bonnet formula for manifolds with boundary [15]; 
it is easily seen that the boundary contribution decays to 0 as the 
boundary is taken to infinity. 

  For an Einstein metric $g$ as in (1.1), $z = 0$, and thus (4.5) gives 
(1.7), via the normalization (1.1). Further, none of the Einstein 
metrics constructed above is conformally flat, i.e. the Weyl tensor $W$ 
does not vanish identically. This is because a conformally flat 
Einstein metric is of constant curvature; however, none of the 
manifolds $M_{\bar \sigma}$ admit a negatively curved metric, as noted 
following (2.5). It follows that for any Dehn filling,
$$vol M_{\bar \sigma} < vol N,$$
see (1.8). Thus, all Einstein manifolds $(M, g)$ obtained by performing 
Dehn filling on the ends of a complete hyperbolic 4-manifold $(N, 
g_{-1})$ have volume less than the volume of $(N, g_{-1})$.

\medskip

 If $N$ is a complete non-compact hyperbolic 4-manifold of finite 
volume, then (4.1) gives
\begin{equation}\label{e4.6}
vol N = \frac{4\pi^{2}}{3}k,
\end{equation}
where $k = \chi(N) \in {\mathbb Z}^{+}$. Thus the volume spectrum of 
hyperbolic 4-manifolds is contained in the set $(4\pi^{2}/3){\mathbb 
Z}^{+}$. 

  Currently, one does not have a complete classification of the 
hyperbolic 4-manifolds of minimal volume $4\pi^{2}/3$, i.e. of Euler 
characteristic 1. However, in [29], an explicit description of 1171 
complete non-compact hyperbolic 4-manifolds is given, all of minimal 
volume $4\pi^{2}/3$. To be concrete, we base the discussion to follow 
on this collection of hyperbolic 4-manifolds, although it is easily 
seen to apply to any initially given hyperbolic 4-manifold.

  Let $N_{a}$, $1 \leq a \leq 1171$ denote the list of complete, 
non-compact hyperbolic 4-manifolds in [29]; of these, 22 are 
orientable, while the rest are non-orientable. Most of the manifolds 
$N_{a}$ have non-zero first Betti number. Hence, for any $k \in 
{\mathbb Z}^{+}$, there are coverings of such manifolds of degree $k$, 
and thus of Euler characteristic $k$ and volume $4\pi^{2}k/3$. It 
follows that the volume spectrum of hyperbolic 4-manifolds is precisely 
the positive integral multiples of $4\pi^{2}/3$. Again the number of 
such distinct manifolds of volume $4\pi^{2}k/3$ grows 
super-exponentially in $k$, as in (1.3). 

\medskip

  All of the manifolds $N_{a}$ above have either 5 or 6 cusp ends. 
However, no $N_{a}$ has all ends given by 3-tori $T^{3}$, (although 
many such $N_{a}$ have double covers with all ends toral). Thus, one 
needs to use Corollary 3.11 to perform Dehn filling on a non-toral end. 
For this, one needs to understand the structure of compact flat 
3-manifolds.

  The classification of compact flat 3-manifolds, cf. [21] or [34] 
shows that there are exactly 10 topological types, 
6 orientable and 4 non-orientable. The 6 orientable manifolds are 
labelled $A$-$F$ in [21] and [29], corresponding to $G_{1}$-$G_{6}$ in 
[34], while the remaining 4 non-orientable manifolds are labelled 
$G$-$J$ in [21], [29] corresponding to $B_{1}$-$B_{4}$ in [34]. The 
3-torus $T^{3}$ corresponds to $A = G_{1}$. Further, the moduli of flat 
structures on such manifolds is completely classified, cf. [34].

  Using the criterion (3.49), a straightforward inspection in [34] 
shows that, among the 10 flat manifolds, only the manifolds $A, B, G, 
H$, (corresponding to $G_{1}, G_{2}, B_{1}, B_{2}$), admit an infinite 
sequence of admissible Dehn fillings. In the notation of [34], $\sigma$ 
may be any primitive (integer coefficient) vector in the plane $\langle 
a_{2}, a_{3} \rangle$ in the case of $G_{2}$, while it may be any such 
vector in the plane $\langle a_{1}, a_{2} \rangle$ in the case of 
$B_{1}$ or $B_{2}$. 

\medskip

   Thus, by Corollary 3.11, infinite sequences of Dehn fillings may be 
applied to any of the cusp ends of the form $A$, $B$, $G$ or $H$, of 
any of the manifolds $N_{a}$, to give complete finite volume Einstein 
metrics. For concreteness, let us illustrate the volume and convergence 
behavior on a specific seed manifold. 

   Take for instance $N_{23}$ from [29]. This manifold has five 
cusp ends, of the type AAGGH, i.e. two of the cusp ends are 3-tori, two 
are of type $G$ and one is of type $H$. The first Betti number of 
$N_{23}$ is given by $b_{1}(N_{23}) = 4$. 

  There are now a number of ways to close off the cusps by Dehn filling.

 (1). Close off any one cusp end of $N_{23}$. This gives an infinite 
sequence of complete Einstein manifolds $(M_{i}^{1}, g_{i}^{1})$, with 
4 cusp ends, converging to $N_{23}$ in the pointed Gromov-Hausdorff 
topology. Formally, to each toral end there are $3\cdot \infty$ 
admissible Dehn fillings, while to each end of type $G$ or $H$, there 
are $1\cdot \infty$ admissible Dehn fillings. With respect to a 
suitable labeling, the volume of $(M_{i}^{1}, g_{i}^{1})$ increases to 
$vol N_{23} = 4\pi^{2}/3$. 
 
 (2). Next, close off any two cusp ends of $N_{23}$, giving a 
(bi)-infinite sequence of complete Einstein manifolds $(M_{i}^{2}, 
g_{i}^{2})$, with 3 cusp ends. If one chooses a subsequence 
$(M_{i'}^{2}, g_{i'}^{2})$ of $(M_{i}^{2}, g_{i}^{2})$ for which both 
Dehn fillings tend to infinity, then $(M_{i'}^{2}, g_{i'}^{2})$ 
converges to $N_{23}$ again and the volumes of $(M_{i'}^{2}, 
g_{i'}^{2})$ increase to $vol N_{23}$, (w.r.t. a suitable labeling). 

 However, if the Dehn filling is fixed on one end, and taken to 
infinity on the other, then the corresponding subsequence of 
$(M_{i}^{2}, g_{i}^{2})$ converges to a complete finite volume Einstein 
metric $(M_{\infty}^{1}, g_{\infty}^{1})$ on a manifold with 4 cusp 
ends. By Theorem 4.3, $(M_{\infty}^{1}, g_{\infty}^{1})$ is one of the 
manifolds constructed in (1) above. While the volumes of the 
subsequence of $(M_{i}^{2}, g_{i}^{2})$ converge to the volume of the 
limit $(M_{\infty}^{1}, g_{\infty}^{1})$, it is not known if this 
convergence can be made monotone, unless the limit $(M_{\infty}^{1}, 
g_{\infty}^{1})$ is hyperbolic, (see below). Since there are infinitely 
many possibilities for the limit $(M_{\infty}^{1}, g_{\infty}^{1})$, 
most of these limits cannot be hyperbolic. 

 (3). Next, close off any three cusp ends of $N_{23}$, giving a 
(tri)-infinite sequence of complete Einstein manifolds $(M_{i}^{3}, 
g_{i}^{3})$, with 2 cusp ends. Limits of sequences in this family are 
then complete, finite volume Einstein manifolds with 3, 4 or 5 ends, of 
the type in (2), (1) or $N_{23}$ respectively.

 (4). Close off any four cusp ends of $N_{23}$, giving a family 
$(M_{i}^{4}, g_{i}^{4})$ with the same features as before.

 (5). Finally, one may close off all 5 cusp ends of $N_{23}$ at once, 
giving a (5-fold)-infinite sequence of compact Einstein manifolds 
$(M_{i}^{5}, g_{i}^{5})$. By taking various subsequences, one obtains 
limits of the form in (1)-(4) above, or again $N_{23}$. 

\medskip

   Thus, one sees that there is a large number of sequences of Einstein 
manifolds, compact or non-compact, converging to the initial seed 
manifold $N_{23}$, as well as many other sequences converging to other 
Einstein limits. The same structure of convergence holds with respect 
to any initial hyperbolic seed manifold $N^{k}$. 

  The discussion above proves the following:
\begin{proposition}\label{p4.4}
Let $N^{k}$ be a complete non-compact hyperbolic 4-manifold, with 
volume $V^{k} = (4\pi^{2}/3)k$, and with $q$ cusps, each of type $A$, 
$B$, $G$, or $H$. Then $N^{k}$ is a $q$-fold limit point of elements of 
${\mathcal E}$, while $V^{k}$ is a $q$-fold limit point of elements in 
${\mathcal V} = vol{\mathcal E} \subset {\mathbb R}^{+}$. 
\end{proposition}
{\endproof}

  An obvious modification of Proposition 4.4 holds when some ends of 
$N^{k}$ are not of the type $A$, $B$, $G$, or $H$.

  Unlike the situation with the Thurston theory in 3 dimensions, it is 
not clear that the volume spectrum ${\mathcal V}$ is well-ordered, (as 
a subset of ${\mathbb R}^{+}$), or finite to one. For the approximate 
metrics $\widetilde g$, although a Dehn-filled end has less volume than 
the corresponding hyperbolic cusp with the same boundary, the difference 
is on the order of $O(R^{-(n-1)})$, which is of the same order as the 
deviation of the Einstein metric $g$ from $\widetilde g$, cf. (3.11). 
Hence, more refined estimates are needed to see if the volume is 
essentially monotone on sequences which open a cusp. 

\bigskip

{\bf \S 4.3.}  Similar results regarding the volume behavior hold at 
least in all even dimensions $n = 2m$. Thus, the Chern-Gauss-Bonnet 
formula in this case states that 
\begin{equation}\label{e4.7}
\chi(M) = \frac{(-1)^{m}}{4^{m}\pi^{m}m!}\int_{M}\sum 
\varepsilon_{i_{1}... i_{n}}R_{i_{1}i_{2}}\wedge ... \wedge 
R_{i_{n-1}i_{n}},
\end{equation}
where the sum is over all permutations of $(1, ..., n)$ and $R$ denotes 
the curvature tensor. This formula holds for all compact manifolds, and 
non-compact manifolds of finite volume of the type considered here. For 
Einstein metrics of the form (1.1), the trace-free part of the Ricci 
curvature vanishes, and $R$ may be written as
\begin{equation}\label{e4.8}
R_{i_{a}i_{b}} = \theta_{i_{a}}\wedge \theta_{i_{b}} + W_{i_{a}i_{b}},
\end{equation}
where $W$ is the Weyl tensor and $\{\theta_{i}\}$ run over an 
orthonormal basis. (Here the sign convention is such that $\langle 
R_{i_{a}i_{b}}, \theta_{i_{b}}\wedge \theta_{i_{a}} \rangle$ gives the 
sectional curvature $K_{i_{a}i_{b}}$). Substituting (4.8) in (4.7) gives
\begin{equation}\label{e4.9}
\chi(M) = \frac{(-1)^{m}2m!}{4^{m}\pi^{m}m!} vol M + \int_{M}P^{m}(W),
\end{equation}
where $P^{m}(W)$ is a polynomial of order $m$ in the Weyl tensor $W$. 
By the same arguments as in \S 4.2, the term $P^{m}(W)$ is small, by 
construction, and becomes arbitrarily small whenever all Dehn fillings 
are sufficiently large. In particular, as the Dehn fillings of each end 
are taken to infinity, one has 
\begin{equation}\label{e4.10}
vol M_{\bar \sigma} \rightarrow vol N = 
(-4\pi)^{m}\frac{m!}{2m!}\chi(N).
\end{equation}

  However, in contrast to the situation in 4-dimensions (1.8), it is 
not known if the term $P^{m}(W)$ has a sign. Hence, it is not known if 
the convergence (4.10) is monotone increasing or decreasing.

  The analogue of Proposition 4.4 regarding the structure of ${\mathcal E}$ 
holds in all dimensions, while the analogue regarding the structure of 
${\mathcal V}$ holds at least in all even dimensions. 

\begin{remark} \label{r4.5}
{\rm An analogue of Theorem 1.1 also holds for complete, conformally 
compact hyperbolic manifolds $(N, g_{-1})$ with a finite number of cusp 
ends, cf. [16].  Such manifolds are of infinite volume, with a finite number 
of expanding ends in addition to the cusp ends. Each expanding end may be 
conformally compactified by a smooth defining function $\rho$ as in 
(2.9). The conformal infinity is then a compact manifold $\partial N$, 
possibly disconnected, with a conformally flat metric $g_{\infty}$. In 
the terminology of Kleinian groups, such manifolds are geometrically 
finite hyperbolic manifolds with a finite number of parabolics.

  Theorem 1.1 generalizes to this context to give the following: a  
sufficiently large Dehn filling of the cusp ends of $(N, g_{-1})$ 
carries a conformally compact Einstein metric $(M, g)$, with the same 
conformal infinity as $(N, g_{-1})$. Consequently, for any such $N$, 
there exist infinitely many conformally compact Einstein manifolds 
$M = M_{\bar \sigma}$, of distinct topological type, which 
have the same conformal infinity $(\partial N, g_{\infty})$. 
We refer to [16] for further details. }
\end{remark}

\section*{Appendix A}
\setcounter{equation}{0}
\begin{appendix}
\setcounter{section}{1}

  In this Appendix, we describe the form of $T^{n-1}$-invariant 
infinitesimal Einstein deformations of the hyperbolic cusp metric; 
this is used to verify the statement (3.28) and computations in Lemma 3.4 
in the proof of Proposition 3.2. 

  Recall that the hyperbolic cusp metric $(C, g_{-1})$ is given by
\begin{equation}\label{eA.1}
g_{C} = r^{-2}dr^{2} + r^{2}g_{T^{n-1}}. 
\end{equation}
An infinitesimal Einstein deformation of $g_{-1}$ is a symmetric 
bilinear form $h$ such that $h \in {\rm Ker}L$, i.e. 
\begin{equation}\label{eA.2}
L(h) = D^{*}Dh - 2R(h) = 0.
\end{equation}
By Lemma 2.1, we need only consider $h$ such that $tr h = 0$. Hence, 
from (2.13) one has
\begin{equation}\label{eA.3}
R(h) = h. 
\end{equation}
Since $h$ is $T^{n-1}$ invariant, $h$ has the form
\begin{equation}\label{eA.4}
h = \sum h_{ij}(r)\theta_{i}\cdot  \theta_{j}, 
\end{equation}
where $\theta_{i}$ is a local orthonormal coframing, dual to $e_{i}$, 
defined as follows: $e_{1} = \nabla s,$ where $ds = r^{-1}dr$, so the 
integral curves of $\nabla s$ are geodesics, while $e_{i}$, $i \geq  2$ 
are tangent to $T^{n-1}$ . If one writes $r^{2}g_{T^{n-1}} = 
r^{2}(d\phi_{2}^{2}+...+d\phi_{n}^{2})$, then $e_{i} = r^{-1}\partial 
/\partial\phi_{i}$ and so $\theta_{i} = rd\phi_{i}$. 

Now we compute $D^{*}Dh = -\nabla_{e_{i}}\nabla_{e_{i}}h + 
\nabla_{\nabla_{e_{i}}e_{i}}h$. From (A.4), one has
$$-\nabla_{e_{i}}\nabla_{e_{i}}h=-\nabla_{e_{i}}\nabla_{e_{i}}
(h_{ab}\theta_{a}\cdot  \theta_{b})= - e_{i}e_{i}(h_{ab})\theta_{a}\cdot  
\theta_{b}-2e_{i}(h_{ab})\nabla_{e_{i}}(\theta_{a}\cdot \theta_{b})- 
h_{ab}\nabla_{e_{i}}\nabla_{e_{i}}(\theta_{a}\cdot \theta_{b}), $$
while
$$\nabla_{\nabla_{e_{i}}e_{i}}h = 
\nabla_{\nabla_{e_{i}}e_{i}}(h_{ab}\theta_{a}\cdot  \theta_{b}) = 
(\nabla_{e_{i}}e_{i})(h_{ab})\cdot (\theta_{a}\cdot \theta_{b}) + 
h_{ab}\nabla_{\nabla_{e_{i}}e_{i}}(\theta_{a}\cdot \theta_{b}). $$
By (A.3), one needs only to consider the $\theta_{a}\cdot  \theta_{b}$ 
component of this. Clearly, by orthogonality of the basis
$$\langle \nabla_{e_{i}}(\theta_{a}\cdot \theta_{b}), \theta_{a}\cdot  
\theta_{b} \rangle = 0 \ {\rm and} \  \langle 
\nabla_{\nabla_{e_{i}}e_{i}}(\theta_{a}\cdot \theta_{b}), 
\theta_{a}\cdot  \theta_{b} \rangle = 0. $$
Combining this with (A.3) and (A.2) then gives
\begin{equation}\label{eA.5}
-\Delta h_{ab} - h_{ab} \langle 
\nabla_{e_{i}}\nabla_{e_{i}}(\theta_{a}\cdot \theta_{b}), 
\theta_{a}\cdot  \theta_{b} \rangle -2h_{ab} = 0. 
\end{equation}

 For $h = h_{ab}$, $h = h(r) = h(s)$, with $dr/ds = r$. Thus
$$\Delta h(s) = (dh/ds)\Delta s + (dh^{2}/ds^{2}). $$
But $dh/ds = h'\cdot (dr/ds) = h' r$, and $dh^{2}/ds^{2} = h' r + h'' 
r^{2}$, with $'  = d/dr$. Also
$$\Delta s = \langle \nabla_{e_{i}}e_{1}, e_{i} \rangle = n-1.$$
Thus,
$$\Delta h(s) = (n-1)rh'  +(rh'  + r^{2}h'' ) = r^{2}h''  + nrh' . $$
Next, one easily computes that:
\begin{equation}\label{eA.6}
\nabla_{e_{1}}\theta_{a} = 0, \ {\rm for \ any} \ a,
\end{equation}
\begin{equation}\label{eA.7}
\nabla_{e_{i}}\theta_{a} = -\delta_{ia}\theta_{1}, \ {\rm for \ any} \ 
a, i > 1, \ {\rm while} \ \nabla_{e_{i}}\theta_{1} = \theta_{i}, \ i > 
1.
\end{equation} 
The latter equations come from fact that the tori are totally umbilic, 
with $2^{\rm nd}$ fundamental form $A = g$, while the intrinsic 
connection on tori is the flat connection, so tangential covariant 
derivatives vanish.

 To compute $\langle \nabla_{e_{i}}\nabla_{e_{i}}(\theta_{a}\cdot 
\theta_{b}), \theta_{a}\cdot  \theta_{b} \rangle$, one has 
$\nabla_{e_{i}}(\theta_{a}\cdot  \theta_{b}) = 
(\nabla_{e_{i}}\theta_{a})\cdot  \theta_{b}+ \theta_{a}\cdot 
\nabla_{e_{i}}\theta_{b}$, and so
$$\nabla_{e_{i}}\nabla_{e_{i}}\theta_{a}\cdot  \theta_{b} = 
(\nabla_{e_{i}}\nabla_{e_{i}}\theta_{a})\cdot  \theta_{b}+ 
2\nabla_{e_{i}}\theta_{a}\cdot \nabla_{e_{i}}\theta_{b} + 
\theta_{a}\cdot \nabla_{e_{i}}\nabla_{e_{i}}\theta_{b}. $$
Suppose first $a > 1$. Then $\nabla_{e_{i}}\theta_{a} = 
-\delta_{ia}\theta_{1},$ so $\nabla_{e_{i}}\nabla_{e_{i}}\theta_{a} = 
-\delta_{ia}\nabla_{e_{i}}\theta_{1} = -\delta_{ia}\theta_{i} = 
-\theta_{a},$ while $\nabla_{e_{i}}\theta_{a} = -\delta_{ia}\theta_{1}
.$ This then gives
$$\langle \nabla_{e_{i}}\nabla_{e_{i}}\theta_{a}\cdot  \theta_{b}, 
\theta_{a}\cdot  \theta_{b} \rangle = -2, \ a, b > 1. $$
Thus, the last two terms in (A.5) cancel and, for $h = h_{ab}$, $a, b > 
1$, one is left with
\begin{equation}\label{eA.8}
\Delta h = 0, \ {\rm i.e.} \  r^{2}h''  + nrh'  = 0.
\end{equation}
The general solution of (A.8) is 
\begin{equation} \label{eA.9}
h = c_{1}r^{-(n-1)} + c_{2}, 
\end{equation}
as in (3.34). 

  Next suppose $a = 1$, $b > 1$, and let $h = h_{1b}$. Then 
$\nabla_{e_{i}}\nabla_{e_{i}}\theta_{1} = -(n-1)\theta_{1}$ and 
$\nabla_{e_{i}}\nabla_{e_{i}}\theta_{b} = -\theta_{b}$. Using 
(A.6)-(A.7) for the middle term in (A.5) then gives $$\langle 
\nabla_{e_{i}}\nabla_{e_{i}}\theta_{1}\cdot  \theta_{b}, 
\theta_{1}\cdot  \theta_{b} \rangle = -(n+2).$$
This gives the Euler equation
\begin{equation}\label{eA.10}
r^{2}h'' + nrh' - nh = 0,
\end{equation}
which has the general solution 
\begin{equation}\label{eA.11}
h = h_{1b} = c_{1}r + c_{2}r^{-n},
\end{equation}
for some constants $c_{1}$, $c_{2}$. 

  Performing similar calculations on $h = h_{11}$ gives the Euler equation 
\begin{equation} \label{eA.12}
r^{2}h + nrh' - 2(n-1)h = 0,
\end{equation}
with general solution 
\begin{equation}\label{eA.13}
h_{11} = c_{1}r^{\alpha_{+}} + c_{2}r^{\alpha_{-}},
\end{equation}
where $\alpha_{\pm} = \frac{1}{2}(-(n-1) \pm \sqrt{(n-1)^2 + 8(n-1)})$.

\end{appendix}

\bibliographystyle{plain}

\bigskip

\begin{center}
March 2003/November, 2005 
\end{center}

\medskip
\noindent
\address
\noindent
{Dept. of Mathematics,\\ 
S.U.N.Y. at Stony Brook\\
Stony Brook, N.Y. 11794-3651\\}
\email{anderson@math.sunysb.edu}

\end{document}